\documentclass[11pt]{article}

\usepackage{color}
\usepackage{latexsym}
\usepackage{dsfont}
\usepackage{amssymb}
\usepackage{graphicx}
\usepackage{amsmath, amsfonts,amssymb,theorem,euscript,array,enumerate,amsfonts,mathrsfs}
\usepackage{hyperref}

\newtheorem{Theorem}{Theorem}[part]
\newtheorem{Definition}{Definition}[part]
\newtheorem{Proposition}{Proposition}[part]

\newtheorem{Lemma}{Lemma}[part]

\def \trans{^{\scriptscriptstyle{\intercal}}}

\def \Int{\displaystyle\int}

\def \b1{\bf{1}}

\def \ba{\bf{a}}

\def \R{\mathbb{R}}
\def \L{\mathbb{L}}
\def \M{\mathbb{M}}

\def \E{\mathbb{E}}
\def \F{\mathbb{F}}

\def \P{\mathbb{P}}
\def \S{\mathbb{S}}

\def \D{\mathbb{D}}

\def\esssup_#1{\underset{#1}{\mathrm{ess\,sup\, }}}
\def\essinf_#1{\underset{#1}{\mathrm{ess\,inf\, }}}

\def\argmin_#1{\underset{#1}{\mathrm{argmin\, }}}

\def \Ac{{\cal A}}

\def \Cc{{\cal C}}
\def \Dc{{\cal D}}

\def \Fc{{\cal F}}
\def \Gc{{\cal G}}

\def \Kc{{\cal K}}
\def \Lc{{\cal L}}
\def \Pc{{\cal P}}

\def \Sc{{\cal S}}

\def \eps{\varepsilon}

\def \ep{\hbox{ }\hfill$\Box$}

\def\reff#1{{\rm(\ref{#1})}}

\def\beqs{\begin{eqnarray*}}
\def\enqs{\end{eqnarray*}}
\def\beq{\begin{eqnarray}}
\def\enq{\end{eqnarray}}

\begin{document}

 \title{Linear quadratic optimal control of conditional McKean-Vlasov equation with random coefficients and applications
 \thanks{This work is part of the ANR project CAESARS (ANR-15-CE05-0024), and also  supported by FiME (Finance for Energy Market Research Centre) and
the ``Finance et D\'eveloppement Durable - Approches Quantitatives'' EDF - CACIB Chair.
 }
 }

\author{Huy\^en PHAM
\\\small  Laboratoire de Probabilit\'es et
 \\\small  Mod\`eles Al\'eatoires, CNRS, UMR 7599
 \\\small  Universit\'e Paris Diderot
 \\\small  pham at math.univ-paris-diderot.fr
\\\small  and CREST-ENSAE
}

\maketitle

\date{}

\begin{abstract}
We consider the optimal control problem for a  linear  conditional McKean-Vlasov equation with quadratic cost functional.  The coefficients of the 
system and the weigh\-ting matrices in the cost functional are allowed to be adapted processes with respect to the common noise filtration. Semi closed-loop strategies are introduced, and following the dynamic programming approach in \cite{phawei16}, we solve the problem and characterize  time-consistent optimal control by means of a system of decoupled backward stochastic Riccati differential equations.  We present several financial applications with explicit solutions, and revisit, in particular, optimal tracking problems with price impact, and the conditional mean-variance portfolio selection in an incomplete market model. 
\end{abstract}

\vspace{5mm}


\vspace{5mm}


\newpage

\section{Introduction and problem formulation}

\setcounter{equation}{0} \setcounter{Assumption}{0}
\setcounter{Theorem}{0} \setcounter{Proposition}{0}
\setcounter{Corollary}{0} \setcounter{Lemma}{0}
\setcounter{Definition}{0} \setcounter{Remark}{0}

Let us formulate the linear quadratic  optimal control of conditional (also called stochastic) 
McKean-Vlasov  equation with random coefficients (LQCMKV in short form).   
Consider the controlled stochastic  McKean-Vlasov dynamics in $\R^d$  given by
\beq 
dX_t &=& b_t(X_t,\E[X_t|W^0],\alpha_t) dt + \sigma_t(X_t,\E[X_t|W^0],\alpha_t) dW_t \nonumber \\
& & \;\;\;\;\; +  \; \sigma^0_t(X_t,\E[X_t|W^0],\alpha_t) dW^0_t, \;\;\; 0 \leq t \leq T,  \;\;\;  X_0 \; = \; \xi_0. \label{stoMcKean}
\enq
Here $W,W^0$ are two independent one-dimensional Brownian motions on some proba\-bility space $(\Omega,\Fc,\P)$,  $\F^0$ $=$ $(\Fc_t^0)_{0\leq t \leq T}$ is the natural filtration generated by $W^0$, 
$\F$ $=$ $(\Fc_t)_{0\leq t\leq T}$ is the natural filtration generated by $(W,W^0)$, augmented with an independent $\sigma$-algebra $\Gc$,  $\xi_0$ $\in$ $L^2(\Gc;\R^d)$ is a square-integrable $\Gc$-measurable random variable with values in 
$\R^d$, $\E[X_t|W^0]$ denotes the conditional expectation of $X_t$ given the whole $\sigma$-algebra $\Fc_T^0$ of $W^0$,  and the control process $\alpha$ is an $\F^0$-progressively measurable process with values in 
$A$ equal either to $\R^m$ or to $L(\R^d;\R^m)$ the set of Lipschitz functions from $\R^d$ into $\R^m$. This distinction of the control sets will be discussed  later  in  the introduction, but for the moment, one may interpret roughly the case when $A$ $=$ $\R^m$ as the modeling for open-loop control and the case when $A$ $=$ $L(\R^d;\R^m)$ as the modeling for closed-loop control. 
When $A$ $=$ $\R^m$, we require that $\alpha$ satisfies the square-integrability condition $L^2(\Omega\times [0,T])$, i.e., 
$\E[\int_0^T |\alpha_t|^2 dt]$ $<$ $\infty$, and we denote by $\Ac$ the set of control processes.  The coefficients $b_t(x,\bar x,a)$, $\sigma_t(x,\bar x,a)$, $\sigma_t^0(x,\bar x,a)$, $0\leq t\leq T$, are $\F^0$-adapted processes with values in $\R^d$, for any $x,\bar x$ $\in$ 
$\R^d$, $a$ $\in$  $A$,  and of linear form:
\begin{equation} \label{XLQ}
\begin{array}{ccc}
b_t(x,\bar x,a) &=& 
\left\{
\begin{array}{ll}
b_t^0 + B_t x + \bar B_t \bar x + C_t a & \mbox{ if } A = \R^m \\
b_t^0 + B_t x + \bar B_t \bar x + C_t a(x)  & \mbox{ if } A =  L(\R^d;\R^m) 
\end{array}
\right. \\
\sigma_t(x,\bar x,a) &=& 
\left\{
\begin{array}{ll}
\gamma_t + D_t x + \bar D_t \bar x + F_t a & \mbox{ if } A = \R^m \\
\gamma_t + D_t x + \bar D_t \bar x + F_t a(x)  & \mbox{ if } A =  L(\R^d;\R^m) 
\end{array}
\right. \\
\sigma_t^0(x,\bar x,a) &=& 
\left\{
\begin{array}{ll}
\gamma_t^0 + D_t^0 x + \bar D_t^0 \bar x + F_t^0 a & \mbox{ if } A = \R^m \\
\gamma_t^0 + D_t^0 x + \bar D_t^0 \bar x + F_t^0 a(x)  & \mbox{ if } A =  L(\R^d;\R^m), 
\end{array}
\right. 
\end{array}
\end{equation}
where $b^0$, $\gamma$, $\gamma^0$ are $\F^0$-adapted processes vector-valued in $\R^d$,  satisfying a square-integrability condition $L^2(\Omega\times [0,T])$: 
$\E[\int_0^T |b_t|^2 + |b_t^0|^2 + |\gamma_t|^2 + |\gamma_t^0|^2 dt]$ $<$  $\infty$,  $B$, $\bar B$, $D$, $\bar D$, $D^0$, $\bar D^0$  
are essen\-tially bounded $\F^0$-adapted processes matrix-valued in $\R^{d\times d}$,  and $C$, $F$, $F^0$ are essen\-tially bounded $\F^0$-adapted processes matrix-valued in $\R^{d\times m}$.  For any $\alpha$ $\in$ $\Ac$,  there exists a unique strong solution $X$ $=$ $X^\alpha$ to \reff{stoMcKean}, which is $\F$-adapted, and satisfies the square-integrability condition $\Sc^2(\Omega\times[0,T])$: 
\beq \label{estimX2}
\E\big[ \sup_{0\leq t\leq T} |X_s^\alpha|^2\big]  & \leq & C_\alpha \big( 1 + \E|\xi_0|^2 \big) \; < \; \infty,  
\enq
for some positive constant $C_\alpha$ depending on $\alpha$: when $A$ $=$ $\R^m$, $C_\alpha$ depends on $\alpha$ via $\E[\int_0^T |\alpha_t|^2 dt]$ $<$ $\infty$, and when $A$ $=$ $L(\R^d;\R^m)$, $C_\alpha$ depends on $\alpha$ via its Lipschitz constant. 

The cost functional to be minimized over $\alpha$ $\in$ $\Ac$ is:
\beqs
J(\alpha) &=& \E \Big[ \int_0^T f_t(X_t^\alpha,\E[X_t^\alpha|W^0],\alpha_t) dt + g(X_T^\alpha,\E[X_T^\alpha|W^0]) \Big], \\
\rightarrow \;\;\; V_0 &  := &  \inf_{\alpha\in\Ac} J(\alpha), 
\enqs 
where $\{f_t(x,\bar x,a)$, $0\leq t\leq T\}$, is an $\F^0$-adapted real-valued process, $g(x,\bar x)$ is a $\Fc_T^0$-measurable random variable, for any $x,\bar x$ $\in$ $\R^d$, $a$ $\in$ $A$, of quadratic form:
\begin{equation} \label{fgquadra}
\begin{array}{ccl} 
f_t(x,\bar x,a) &=& \left\{
\begin{array}{ll}
x\trans Q_t x +  \bar x\trans \bar Q_t \bar x + M_t\trans x + a\trans N_t a   & \mbox{ if } A = \R^m \\
x\trans Q_t x +  \bar x\trans \bar Q_t \bar x +  M_t\trans x +  a(x)\trans N_t a(x)    & \mbox{ if } A =  L(\R^d;\R^m) 
\end{array}
\right.   \\
g(x,\bar x) &=&  x\trans P x +  \bar x\trans \bar P \bar x + L\trans x,   
\end{array}
\end{equation}
where $Q$, $\bar Q$ are essentially bounded $\F^0$-adapted processes, with values in $\S^d$ the set of symmetric matrices in $\R^{d\times d}$,   
$P$, $\bar P$ are essentially bounded $\Fc_T^0$-measurable random matrices in $\S^d$, $N$  is an essentially bounded $\F^0$-adapted process, with values in $\S^m$, 
$M$ is an  $\F^0$-adapted process with values in $\R^d$, satisfying a square integrability condition $L^2(\Omega\times [0,T])$,  $L$ is an $\Fc_T^0$-measurable square integrable random vector in $\R^d$, 
and $\trans$ denotes the transpose of any vector or matrix.   



\vspace{1mm}

The above control formulation  of  stochastic  McKean-Vlasov equations provides  a unified framework for  some important classes of control problems. 
In particular, it is motivated in particular  by the asymptotic formulation of cooperative equilibrium for a large population of particles (players) in  mean-field 
interaction under common noise (see, e.g., \cite{cardel13}, \cite{carzhu14}) and also occurs when the cost functional involves the first and second moment of the (conditional) law of the state process, for example in (conditional) mean-variance portfolio selection problem (see, e.g., \cite{lizho00}, \cite{bascha10}, \cite{borkum10}).  
When $A$ $=$ $L(\R^d;\R^m)$,  this corresponds to the problem of a (representative) agent, using a control $\alpha$ based on her/his current private state $X_t$ at time $t$, and of the information brought by the common noise $\Fc_t^0$, typically the conditional mean $\E[X_t|W^0]$, which represents,   
in the large population equilibrium interpretation, the  limit   of the em\-pirical mean of the state of all the players when their number tend to infinity from  the propagation of chaos. 
In other words, the control $\alpha$ may be viewed as a  semi closed-loop control, i.e., closed-sloop w.r.t. the state process, and open-loop w.r.t. the common noise $W^0$, or alternatively as a $\F^0$-progressively measurable random field control 
$\alpha$ $=$ $\{\alpha_t(x), 0 \leq t \leq T, x \in \R^d\}$.  This class of semi closed-loop control extends the class of closed-loop strategies for the LQ  control of McKean-Vlasov equations (or mean-field stochastic differential equations) without common noise $W^0$,  as recently studied in  \cite{lietal16} where the controls are chosen at any time $t$ in linear form w.r.t. the current state value $X_t$ and the deterministic expected  value $\E[X_t]$. When $A$ $=$ $\R^m$,  the LQCMKV pro\-blem may be viewed as a special  partial observation control  problem for a state dynamics like in \reff{stoMcKean} where the controls are of open-loop form, and adapted w.r.t. an observation filtration $\F^I$ $=$ $\F^0$ generated by some exogenous random factor process $I$ driven by $W^0$. In the case where $\sigma$ $=$ $0$, we see that  the process $X$ is $\F^0$-adapted,  hence $\E[X_t|W^0]$ $=$ $X_t$, and the LQCMKV  problem is reduced to the classical LQ  control problem (see, e.g., \cite{yonzho99})  with random coefficients,  with  open-loop controls for $A$ $=$ $\R^m$ or closed-loop controls for $A$ $=$ $L(\R^d;\R^m)$.  Note that this distinction between  open-loop and closed-loop strategies for LQ control problems has been recently introduced in \cite{sunyon14} where closed-loop controls are assumed of  linear form w.r.t. the current  state value, while it is considered here  a priori only Lipschitz w.r.t. the current state value.  

Optimal control of McKean-Vlasov equation is a rather new topic in the area of stochastic control and applied probability, and addressed, e.g., in \cite{anddje10}, \cite{bucetal11}, \cite{benetal13},  \cite{cardel14},  \cite{phawei15}. In this McKean-Vlasov context, the class of linear quadratic optimal control,  
which provides a typical case for solvable applications, has been studied in several papers, among them \cite{huetal11},  \cite{yon13},  \cite{huaetal15}, \cite{sun}, where the coefficients are assumed to be deterministic.  
It is often argued that due to the pre\-sence of the  law of the state in a nonlinear way 
(here for the LQ problem, the square of the  expectation),  the problem is time-inconsistent in the sense that an optimal control viewed from today is no more optimal when viewed from tomorrow, and this would prevent a priori the use of the dynamic programming method. To tackle time inconsistency, one then focuses typically on either pre-commitment strategies, i.e. controls that are optimal  for the problem viewed at  the initial time, but may be not optimal at future date, 
or  game-equilibrium strategies, i.e., control decisions  considered as a game against all the future decisions the controller  is going to make.

In this paper,  we shall focus on the optimal control for the initial value $V_0$ of the LQCMKV problem with random coefficients,  but  following the approach developed in \cite{phawei16}, we emphasize that time consistency can be actually restored for pre-commitment strategies, provided that one considers as state variable  the conditional law of the state process instead of the state itself, therefore making possible the use of the dynamic programming method. We show  that the dynamic version of the LQCMKV control problem defined by a random field value function,  has a quadratic structure with respect to the conditional law of the state process, leading to a characterization of the optimal control in terms of a decoupled system of backward stochastic Riccati equations (BSREs) whose existence and uniqueness are obtained in connection with a standard LQ control problem.  
The main ingredient for such derivation is an 
It\^o's formula along a flow of conditional measures and a suitable notion of differentiability with respect to probability measures.  We illustrate our results with several financial applications. We first revisit the optimal trading and  benchmark tracking problem with price impact for  general price and target processes, and obtain closed-form solutions extending some known results in the literature.  We next solve a variation of the mean-variance portfolio selection problem 
in an incomplete market with random factor.  Our last example considers an  interbank systemic risk model with random factor in a common noise environment.

The paper is organized as follows.  Section 2 gives some key preliminaries: we reformulate the LQCMKV problem into a problem  involving the conditional law of the state process as state variable for which a dynamic programming verification theorem is stated and time consistency holds. We also recall the It\^o's formula along a flow of conditional measures.  Section 3 is devoted to the characterization of the optimal control by means of a system of BSREs in the case of both a control set $A$ $=$ $\R^m$ and $A$ $=$ $L(\R^d;\R^m)$.  We develop  in Section 4 the applications.

\vspace{1mm}

We end this introduction with some notations. 

\vspace{1mm}

\noindent {\bf Notations}.  We denote by $\Pc_{_2}(\R^d)$  the set probability measures  $\mu$ on $\R^d$, which are square integrable, i.e., $\|\mu\|_{_2}^2$ $:=$ $\int_{\R^d} |x|^2 \mu(dx)$ $<$ $\infty$. 
For any $\mu$ $\in$ $\Pc_{_2}(\R^d)$,   we denote by $L_\mu^2(\R^q)$ the set of measurable functions $\varphi$ $:$ $\R^d$ $\rightarrow$ $\R^q$, which are square integrable with respect to $\mu$, by $L_{\mu\otimes\mu}^2(\R^q)$ 
the set of measurable functions $\psi$ $:$  $\R^d\times\R^d$ $\rightarrow$ $\R^q$,  which are square integrable with respect to the product measure  $\mu\otimes\mu$, and we set
\beqs
\mu(\varphi) \; := \;   \int \varphi(x)\,\mu(dx),  \;\; \bar\mu \; := \; \int x \mu(dx), & & \mu\otimes\mu(\psi) \;  := \;  \int  \psi(x,x') \mu(dx)\mu(dx'). 
\enqs 
We also define $L_\mu^\infty(\R^q)$ (resp.  $L_{\mu\otimes\mu}^\infty(\R^q)$) as the subset of elements $\varphi$ $\in$ $L_\mu^2(\R^q)$ (resp.  $L_{\mu\otimes\mu}^2(\R^q)$) which are bounded 
$\mu$ (resp. $\mu\otimes\mu$) a.e., and $\|\varphi\|_\infty$ is their essential supremum.  For any random variable $X$ on $(\Omega,\Fc,\P)$, we denote by $\Lc(X)$  its probability law (or distribution) under  $\P$, by 
$\Lc(X|W^0)$ its conditional law given $\Fc_T^0$,   and we shall assume w.l.o.g. that $\Gc$ is rich enough in the sense that  $\Pc_{_2}(\R^d)$ $=$ $\{\Lc(\xi): \xi \in L^2(\Gc;\R^d)\}$.

\section{Preliminaries} \label{secpreli}

\setcounter{equation}{0} \setcounter{Assumption}{0}
\setcounter{Theorem}{0} \setcounter{Proposition}{0}
\setcounter{Corollary}{0} \setcounter{Lemma}{0}
\setcounter{Definition}{0} \setcounter{Remark}{0}

For any $\alpha$ $\in$ $\Ac$, and $X^\alpha$ $=$ $(X_t^\alpha)_{0\leq t\leq T}$ the solution to \reff{stoMcKean}, we define $\rho_t^\alpha$ $=$ $\Lc(X_t^\alpha|W^0)$ as the 
conditional law of $X_t^\alpha$ given $\Fc_T^0$ for $0\leq t\leq T$. Since $X^\alpha$ is $\F$-adapted, and $W^0$ is a $(\P,\F)$-Wiener process, we notice that 
$\rho_t^\alpha(dx)$ $=$ $\P[X_t^\alpha \in dx| \Fc_T^0]$ $=$ $\P[X_t^\alpha \in dx| \Fc_t^0]$, and thus  $\{\rho_t^\alpha,0\leq t\leq T\}$ admits an $\F^0$-progressively measurable modification (see, e.g., Theorem 2.24 in \cite{baicri09}), that will be identified with itself in the sequel, and is valued  in $\Pc_{_2}(\R^d)$ by \reff{estimX2}, namely:
\beq  \label{estimrho2}
\E \big[ \sup_{0\leq t\leq T} \|\rho_t^\alpha\|_{_2}^2 \big] & \leq & C_\alpha \big( 1 + \E|\xi_0|^2 \big).
\enq
Moreover, we mention that the process $\rho^\alpha$ $=$ $(\rho_t^\alpha)_{0\leq t\leq T}$ has continuous trajectories as it is valued in 
$\Pc_{_2}(C([0,T];\R^d)$ the set of square integrable probability measures on the space $C([0,T];\R^d)$ of continuous functions from $[0,T]$ into $\R^d$.

Now, by the  law of iterated conditional expectations, and recalling that $\alpha$ $\in$ $\Ac$ is $\F^0$-progressively measurable, 
we can rewrite the cost functional as
\beq 
J(\alpha) &=& \E \Big[  \int_0^T \E \big[  f_t\big(X_t^{\alpha},\bar\rho_t^\alpha,\alpha_t  \big) \big| \Fc_t^0 \big] dt
+ \E \big[ g\big(X_T^{\alpha},\bar\rho_T^\alpha\big) \big| \Fc_T^0 \big]  \Big] \nonumber \\
&=& \E \Big[ \int_0^T \rho_t^{\alpha} \big( f_t(.,\bar\rho_t^{\alpha},\alpha_t)\big) dt + \rho_T^{\alpha}\big(g(.,\bar\rho_T^{\alpha})\big) \Big] \nonumber \\
&=& \E \Big[ \int_0^T \hat f_t(\rho_t^\alpha,\alpha_t) dt + \hat g(\rho_T^{\alpha}) \Big], \label{Jmu}
\enq
where we used in the second equality the fact that $\{f_t(x,\bar x,a),  x,\bar x \in \R^d, a \in A, 0\leq t\leq T\}$,  is a random field $\F^0$-adapted process, $g(x)$ is $\Fc_T^0$-measurable, 
and the $\F^0$-adapted process  $\{\hat f_t(\mu,a),0\leq t\leq T\}$, the $\Fc_T^0$-measurable random variable  $\hat g(\mu)$, for $\mu$ $\in$ $\Pc_{_2}(\R^d)$, $a$ $\in$ $A$,  are defined by
\begin{equation*}\label{hatfg}
\left\{
\begin{array}{ccc}
\hat f_t(\mu,a) &:=& \mu\big( f_t(.,\bar\mu,a)\big) \; = \; \int f_t(x,\bar \mu,a) \mu(dx) \\
\hat g(\mu) &:=& \mu\big( g(.,\bar\mu) \big) \; = \; \int g(x,\bar \mu) \mu(dx).
\end{array}
\right. 
\end{equation*}
From the quadratic forms of $f,g$ in \reff{fgquadra}, the random fields $\hat f_t(\mu,a)$ and $\hat g(\mu)$, 
$(t,\mu,a)$ $\in$ $[0,T]\times\Pc_{_2}(\R^d)\times A$, are given  by
\begin{equation} \label{hatfgquadra}
\begin{array}{ccll}
\hat f_t(\mu,a) &=&  \left\{
\begin{array}{ll}
{\rm Var}(\mu,Q_t)  + v_2(\mu,Q_t + \bar Q_t) & \\
\; + \;  v_1(\mu,M_t) + a\trans N_t a   & \mbox{ if } A = \R^m \\
{\rm Var}(\mu,Q_t)  + v_2(\mu,Q_t + \bar Q_t) & \\
\; + \;   v_1(\mu,M_t) +   \Int [a(x)\trans N_t a(x)  ] \mu(dx)   & \mbox{ if } A =  L(\R^d;\R^m), \\
\end{array}
\right. \\
\hat g(\mu) &=& {\rm Var}(\mu,P) + v_2(\mu,P+\bar P) + v_1(\mu,L), & 
\end{array}
\end{equation}
where we define the functions on  $\Pc_{_2}(\R^d)\times\S^d$  and   $\Pc_{_2}(\R^d)\times\R^d$ by:
\beqs
{\rm Var}(\mu,k) &:=& \int (x-\bar\mu)\trans k (x-\bar\mu) \mu(dx), \;\;\; \mu \in \Pc_{_2}(\R^d), \; k \in \S^d,      \\
v_2(\mu,\ell) &:=&   \bar\mu\trans \ell \bar\mu, \;\;\; \mu \in \Pc_{_2}(\R^d), \; \ell \in \S^d \\
v_1(\mu,y) &:=&  y\trans\bar\mu,  \;\;\;\;\;\;\;  \mu \in \Pc_{_2}(\R^d), \; y \in \R^d. 
\enqs

We shall make the following assumptions  on the  coefficients of the model:

\vspace{2mm}

\noindent {\bf (H1)} \hspace{1cm} $Q$, $Q+\bar Q$, $P$, $P+\bar P$, $N$ are nonnegative a.s.; 

\vspace{1mm}

\noindent {\bf (H2)}  \hspace{1cm} One of the two following conditions holds: 
\begin{itemize}
\item[(i)] $N$ is uniformly positive definite  i.e. $N_t$ $\geq$ $\delta I_m$, $0\leq t\leq T$, a.s. for some $\delta$ $>$  $0$; 
\item[(ii)] $P$ or $Q$ is uniformly positive definite, and $F$ is uniformly nondegenerate, i.e. $|F_t|$ $\geq$ $\delta$, $0\leq t\leq T$, a.s., for some $\delta$ $>$ $0$.  
\end{itemize}

\vspace{2mm}

Let us define the dynamic formulation of the stochastic McKean-Vlasov control problem. For any $t$ $\in$ $[0,T]$, $\xi$ $\in$ $L^2(\Gc;\R^d)$, and $\alpha$ $\in$ $\Ac$,  there exists a unique strong solution, denoted by 
$\{X_s^{t,\xi,\alpha},t\leq s\leq T\}$,  to the equation \reff{stoMcKean} starting from  $\xi$ at time $t$, and by noting that $X^{t,\xi,\alpha}$ is also unique in law, we see that the conditional law of $X_s^{t,\xi,\alpha}$ given $\Fc_T^0$ 
depends on $\xi$ only through  its law $\Lc(\xi)$ $=$ $\Lc(\xi|W^0)$ (recall that $\Gc$ is independent of $W^0$).  Then, recalling also that  $\Gc$ is rich enough, the relation 
\beqs
\rho_s^{t,\mu,\alpha} &:=& \Lc(X_s^{t,\xi,\alpha} | W^0), \;\;\; t \leq s \leq T, \; \mu = \Lc(\xi), 
\enqs
defines for any $t$ $\in$ $[0,T]$, $\mu$ $\in$ $\Pc_{_2}(\R^d)$, and $\alpha$ $\in$ $\Ac$, an $\F^0$-progressively measurable continuous process (up to a modification) $\{\rho_s^{t,\mu,\alpha},t\leq s\leq T\}$, with values in $\Pc_{_2}(\R^d)$, and 
as a consequence of the pathwise uniqueness of the solution $\{X_s^{t,\xi,\alpha},t\leq s\leq T\}$, we have the flow property for the conditional law (see Lemma 3.1 in \cite{phawei16} for details):
\beq \label{flowrho}
\rho_s^\alpha &=& \rho_s^{t,\rho_t^\alpha,\alpha}, \;\;\; t \leq s \leq T, \; \alpha \in \Ac. 
\enq 
We then consider the conditional cost functional
\beqs
J_t(\mu,\alpha) & = & \E\Big[ \int_t^T \hat f_s(\rho_s^{t,\mu,\alpha},\alpha_s) ds + \hat g(\rho_T^{t,\mu,\alpha}) \big| \Fc_t^0 \Big], \;\;\; t \in [0,T],  \mu \in \Pc_{_2}(\R^d), \alpha\in \Ac, 
\enqs
which is well-defined by \reff{estimrho2} and under the boundedness assumptions on the  weighting matrices of the quadratic cost function. We next define  
the $\F^0$-adapted random field value function
\beqs
v_t(\mu) &=& \essinf_{\alpha\in\Ac} J_t(\mu,\alpha),  \;\;\;  t \in [0,T],  \mu \in \Pc_{_2}(\R^d),
\enqs
so that
\beq \label{V0v}
V_0 & := & \inf_{\alpha\in\Ac} J(\alpha) \; = \;  v_0(\Lc(\xi_0)),
\enq
which may take a priori for the moment the value $-\infty$.  We shall see later that the Assumptions {\bf (H1)} and {\bf (H2)} will ensure that $V_0$ is finite and there exists an optimal control. The dynamic counterpart of \reff{V0v} is given by 
\beq \label{Vtv}
V_t^\alpha &:=& \essinf_{\beta\in\Ac_t(\alpha)} J_t(\rho_t^\alpha,\beta) \; = \; v_t(\rho_t^\alpha), \;\;\;  t \in [0,T], \; \alpha \in \Ac, 
\enq
where $\Ac_t(\alpha)$ $=$ $\{ \beta  \in \Ac: \beta_s = \alpha_s, s \leq t\}$, and the second equality in \reff{Vtv} follows from the flow property \reff{flowrho} and the observation that $\rho_t^\beta$ $=$ $\rho_t^\alpha$ for $\beta$ $\in$ 
$\Ac_t(\alpha)$.  


\vspace{1mm}

By using general results in \cite{elk81} for dynamic programming, one can   show   (under the condition that the random field $v(\mu)$ is finite) that    
the process $\{v_t(\rho_t^\alpha) + \int_0^t \hat f_s(\rho_s^{\alpha},\alpha_s) ds,0\leq t\leq T\}$ 
is a $(\P,\F^0)$-submartingale, for any $\alpha$ $\in$ $\Ac$, and  
$\alpha^*$ $\in$ $\Ac$ is an optimal control for $V_0$ if and only if $\{v_t(\rho_t^{\alpha^*}) + \int_0^t \hat f_s(\rho_s^{\alpha^*},\alpha_s^*) ds,0\leq t\leq T\}$ is a $(\P,\F^0)$-martingale.  
We shall use  a converse result, namely  a dynamic programming verification theorem, which takes the following formulation in our context.

\begin{Lemma} \label{lemverif}
Suppose that one can find an $\F^0$-adapted random field  $\{w_t(\mu), 0\leq t\leq T, \mu\in\Pc_{_2}(\R^d)\}$ satisfying the quadratic growth condition
\beq \label{wquadragrowth}
|w_t(\mu)| & \leq & C\|\mu\|_{_2}^2  + I_t, \;\;\; \mu \in  \Pc_{_2}(\R^d), \; 0 \leq t \leq T, \; a.s. 
\enq
for some positive constant $C$, and nonnegative $\F^0$-adapted process $I$ with  $\E[\sup_{0\leq t\leq T}|I_t|]$ $<$ $\infty$, such that 
\begin{itemize}
\item[(i)] $w_T(\mu)$ $=$ $\hat g(\mu)$, $\mu$ $\in$ $\Pc_{_2}(\R^d)$;
\item[(ii)] $\{w_t(\rho_t^\alpha)+ \int_0^t \hat f_s(\rho_s^{\alpha},\alpha_s) ds ,0\leq t\leq T\}$ is a $(\P,\F^0)$ local submartingale, for any 
$\alpha$ $\in$ $\Ac$;
\item[(iii)] there exists $\hat\alpha$ $\in$ $\Ac$ such that $\{w_t(\rho_t^{\hat\alpha}) + \int_0^t \hat f_s(\rho_s^{\hat\alpha},\hat\alpha_s) ds,0\leq t\leq T\}$ is a $(\P,\F^0)$ local martingale.
\end{itemize}
Then $\hat\alpha$ is an optimal control for $V_0$, i.e. $V_0$ $=$ $J(\hat\alpha)$,  and 
\beqs
V_0 &=&   w_0(\Lc(\xi_0)). 
\enqs
Moreover, $\hat\alpha$ is time consistent in the sense that
\beqs
V_t^{\hat\alpha} &=& J_t(\rho_t^{\hat\alpha},\hat\alpha), \;\;\;  \forall 0 \leq t\leq T. 
\enqs
\end{Lemma}
{\bf Proof.} By  the local submartingale property in condition (ii), there exists a nondecreasing sequence of $\F^0$-stopping times $(\tau_n)_n$, $\tau_n$ $\nearrow$ $T$ a.s., such that
\beq \label{localsub}
\E[w_{\tau_n}(\rho_{\tau_n}^\alpha) + \int_0^{\tau_n}  \hat f_t(\rho_t^\alpha,\alpha_t) dt]   & \geq &   w_0(\rho_0^\alpha) \; = \; w_0(\Lc(\xi_0)), \;\;\; \forall \alpha \in \Ac  
\enq
From the quadratic form of $f$ in \reff{fgquadra}, we easily see that for all $n$, 
\beqs
\E \Big[ \big| \int_0^{\tau_n} \hat f_t(\rho_t^\alpha,\alpha_t) dt \big| \Big] & \leq &  C_\alpha \Big( 1 + \E\big[  \sup_{0\leq t\leq T} \|\rho_t^\alpha\|_{_2}^2 \big]  \Big),
\enqs
for some positive constant $C_\alpha$ depending on $\alpha$ (when $A$ $=$ $\R^m$, $C_\alpha$ depends on $\alpha$ via $\E[\int_0^T |\alpha_t|^2 dt]$ $<$ $\infty$, and when $A$ $=$ $L(\R^d;\R^m)$, $C_\alpha$ depends on $\alpha$ via its Lipschitz constant). Together with the quadratic growth condition of $w$, and from \reff{estimrho2}, one can then apply dominated convergence theorem  by sending $n$ to infinity into \reff{localsub}, and get
\beqs
w_0(\Lc(\xi_0)) & \leq & \E[w_{T}(\rho_{T}^\alpha) + \int_0^{T}  \hat f_t(\rho_t^\alpha,\alpha_t) dt] \; = \; \E[\hat g(\rho_{T}^\alpha) + \int_0^{T}  \hat f_t(\rho_t^\alpha,\alpha_t) dt] \; = \; J(\alpha)
\enqs 
where we used the terminal condition (i), and the expression \reff{Jmu} of the cost functional. Since $\alpha$ is arbitrary in $\Ac$, this shows that $w_0(\Lc(\xi_0))$ $\leq$ $V_0$.  
The equality is obtained with the local martingale property for $\hat\alpha$ in condition (iii).

From the flow property \reff{flowrho}, and since $\rho_t^\beta$ $=$ $\rho_t^{\hat\alpha}$ for $\beta$ $\in$ $\Ac_t(\hat\alpha)$, we notice that the local submartingale and martingale properties in (ii) and (iii) are formulated on the interval 
$[t,T]$ as:   
\begin{itemize}
\item $\{w_s(\rho_s^{t,\rho_t^{\hat\alpha},\beta})+ \int_t^s \hat f_u(\rho_u^{t,\rho_t^{\hat\alpha},\beta},\beta_u) du ,t\leq s\leq T\}$ is a $(\P,\F^0)$ local submartingale, for any 
$\beta$ $\in$ $\Ac_t(\hat\alpha)$;
\item $\{w_s(\rho_s^{t,\rho_t^{\hat\alpha},\hat\alpha}) + \int_t^s \hat f_u(\rho_u^{t,\rho_t^{\hat\alpha},\hat\alpha},\hat\alpha_u) du, t\leq s\leq T\}$ is a $(\P,\F^0)$ local martingale.
\end{itemize}
By the same arguments as for the initial date, this implies that $V_t^{\hat\alpha}$ $=$ $J_t(\rho_t^{\hat\alpha},\hat\alpha)$ $=$ $w_t(\rho_t^{\hat\alpha})$, which means that $\hat\alpha$ is an optimal control over $[t,T]$, once we start at time $t$ from the initial state $\rho_t^{\hat\alpha}$, i.e.,  the time consistency of $\hat\alpha$.  
\ep

\vspace{2mm}

The practical application of Lemma \ref{lemverif} consists in finding a random field $\{w_t(\mu), \mu \in \Pc_{_2}(\R^d),0\leq t\leq T\}$, smooth (in a sense to be precised), so that one can apply an It\^o's formula to  
$\{w_t(\rho_t^\alpha)+ \int_0^t \hat f_s(\rho_s^{\alpha},\alpha_s) ds ,0\leq t\leq T\}$, and check that the finite variation term is nonnegative for any 
$\alpha$ $\in$ $\Ac$ (the local submartingale condition), and equal to zero for some 
$\hat\alpha$ $\in$ $\Ac$ (the local martingale condition).  For this purpose, we need a notion of derivative with respect to a probability measure, and shall rely on the one introduced by P.L. Lions in his course at Coll\`ege de France \cite{lio12}. We briefly recall the basic definitions and refer to  \cite{car12} for the details,  see also \cite{buetal14}, \cite{chacridel15}. 
This notion is based on the {\it lifting}  of functions $u$ defined on  $\Pc_{_2}(\R^d)$ into functions $U$ defined on  $L^2(\Gc;\R^d)$  by setting $U(X)$ $=$  $u(\Lc(X))$.    
We say that $u$ is differentiable (resp. $\Cc^1$) on $\Pc_{_2}(\R^d)$ if the lift $U$ is Fr\'echet differentiable (resp. Fr\'echet differentiable with continuous derivatives) on  $L^2(\Gc;\R^d)$. 
In this case, the Fr\'echet derivative viewed as an element $DU(X)$ of $L^2(\Gc;\R^d)$  by Riesz's theorem  can be represented as
\beqs \label{Uu1}
DU(X) &=& \partial_\mu u(\Lc(X))(X),
\enqs
for some function  $\partial_\mu u(\Lc(X))$ $:$ $\R^d$ $\rightarrow$ $\R^d$,  which is  called derivative of $u$ at $\mu$ $=$ $\Lc(X)$.  Moreover, 
$\partial_\mu u(\mu)$ $\in$ $L^2_\mu(\R^d)$ for  $\mu$ $\in$ $\Pc_{_2}(\R^d)$ $=$ $\{ \Lc(X): X \in L^2(\Gc;\R^d)\}$. 
Following \cite{chacridel15}, we say that $u$ is fully  $\Cc^2$ if it is $\Cc^1$, the mapping 
$(\mu,x)$ $\in$ $\Pc_{_2}(\R^d)\times\R^d$ $\mapsto$ $\partial_\mu u(\mu)(x)$  is continuous and
\begin{itemize}
\item[(i)] for each fixed $\mu$ $\in$ $\Pc_{_2}(\R^d)$, the mapping $x$ $\in$ $\R^d$ $\mapsto$  $\partial_\mu u(\mu)(x)$ is differentiable in the standard sense, with a gradient denoted by  
$\partial_x  \partial_\mu u(\mu)(x)$  $\in$ $\R^{d\times d}$, and s.t. the mapping  $(\mu,x)$ $\in$ $\Pc_{_2}(\R^d)\times\R^d$ 
$\mapsto$  $\partial_x  \partial_\mu u(\mu)(x)$ is continuous;
\item[(ii)] for each fixed  $x$ $\in$ $\R^d$, the mapping $\mu$ $\in$ $\Pc_{_2}(\R^d)$ $\mapsto$  $\partial_\mu u(\mu)(x)$ is differentiable in the above lifted sense.  Its derivative, interpreted thus as a mapping $x'$ $\in$ $\R^d$ $\mapsto$ $\partial_\mu \big[ \partial_\mu u(\mu)(x)\big](x')$ $\in$ $\R^{d\times d}$ in 
$L^2_\mu(\R^{d\times d})$, is denoted by $x'$ $\in$ $\R^d$ $\mapsto$ $\partial_\mu^2 u(\mu)(x,x')$, and s.t. the mapping $(\mu,x,x')$ $\in$ 
$\Pc_{_2}(\R^d)\times\R^d\times\R^d$ $\mapsto$ $\partial_\mu^2 u(\mu)(x,x')$ is continuous. 
\end{itemize}
We say that $u$ $\in$ $\Cc^2_b(\Pc_{_2}(\R^d))$ if it is fully $\Cc^2$,  $\partial_x  \partial_\mu u(\mu)$ $\in$ $L_\mu^\infty(\R^{d\times d})$, $\partial_\mu^2 u(\mu)$ $\in$ $L_{\mu\otimes\mu}^\infty(\R^{d\times d})$ for any $\mu$ 
$\in$ $\Pc_{_2}(\R^d)$,  and for any compact set $\Kc$ of $\Pc_{_2}(\R^d)$, we have
\beqs
 \sup_{ \mu \in \Kc } \Big[ \int_{\R^d} \big| \partial_\mu u(\mu)(x) |^2\mu(dx)  +
\big \| \partial_x \partial_\mu u(\mu)\|_{_\infty} +  \big \| \partial_\mu^2 u(\mu)\|_{_\infty}
\Big]  & < & \infty.
\enqs
We next need an It\^o's formula along a flow of conditional measures proved in \cite{cardel14b} for processes with common noise. In our context, for the flow of the conditional law $\rho_t^\alpha$, $0\leq t\leq T$, $\alpha$ $\in$ $\Ac$,  
it is formulated as follows. Let $u$ $\in$ $\Cc^2_b(\Pc_{_2}(\R^d))$. Then, for all $t$ $\in$ $[0,T]$, we have
\beq
u(\rho_t^\alpha) &=& u(\Lc(\xi_0))  + \int_0^t  \rho_t^\alpha \big( \L_t^{\alpha_t} u(\rho_t^\alpha) \big) + \rho_t^\alpha\otimes\rho_t^\alpha \big( \M_t^{\alpha_t} u(\rho_t^\alpha) \big) dt \nonumber \\
& & \hspace{3cm}  \;\;\; + \; \int_0^t  \rho_t^\alpha\big( \D_t^{\alpha_t} u(\rho_t^\alpha)  \big) dW_t^0, \label{Ito}
\enq
where for $(t,\mu,a)$ $\in$ $[0,T]\times\Pc_{_2}(\R^d)\times A$, $\L_t^a u(\mu)$, $\D_t^a u(\mu)$  are the $\Fc_t^0$-measurable random functions in $L_\mu^2(\R)$ defined by
\beqs
\L_t^a u(\mu)(x) &:=& b_t(x,\bar\mu,a).\partial_\mu u(\mu)(x) + \frac{1}{2}{\rm tr}\big(\partial_x\partial_\mu u(\mu)(x)(\sigma_t\sigma_t\trans + \sigma_t^0(\sigma_t^0)\trans)(x,\bar\mu,a) \big), \\
\D_t^a u(\mu)(x) &:=& \partial_\mu u(\mu)(x)\trans\sigma_t^0(x,\bar\mu,a),
\enqs
and $\M_t^a u(\mu)$ is the $\Fc_t^0$-measurable random function in $L_{\mu\otimes\mu}^2(\R)$ defined by
\beqs
\M_t^a u(\mu)(x,x') &:=& \frac{1}{2}{\rm tr}\big(\partial^2_\mu u(\mu)(x,x')\sigma_t^0(x,\bar\mu,a)(\sigma_t^0)\trans(x',\bar\mu,a) \big).
\enqs     
  
The dynamic programming verification result in Lemma \ref{lemverif} and It\^o's formula \reff{Ito} are valid for a general stochastic McKean-Vlasov 
equation (beyond the LQ framework), and by combining with an It\^o-Kunita type formula  for random field processes, similar to the one in 
\cite{kun82}, one could apply it to $\{w_t(\rho_t^\alpha)+ \int_0^t \hat f_s(\rho_s^{\alpha},\alpha_s) ds ,0\leq t\leq T\}$ in order to derive a form of 
stochastic Hamilton-Jacobi-Bellman, i.e., a backward stochastic partial differential equation (BSPDE) for $w_t(\mu)$, as done in  
\cite{pen92} for controlled diffusion processes with random coefficients. 
We postpone this general approach for further study and, in the next sections, return to the important special case of LQCMKV problem[s] for which we show that BSPDE[s] are reduced to backward stochastic Riccati equations (BSRE) as in the classical LQ framework.


\section{Backward stochastic Riccati equations}

\setcounter{equation}{0} \setcounter{Assumption}{0}
\setcounter{Theorem}{0} \setcounter{Proposition}{0}
\setcounter{Corollary}{0} \setcounter{Lemma}{0}
\setcounter{Definition}{0} \setcounter{Remark}{0}

We search for an $\F^0$-adapted random field solution to the LQCMKV problem in the quadratic form
\beq \label{wLQ}
w_t(\mu) &=& {\rm Var}(\mu,K_t) +  v_2(\mu,\Lambda_t) + v_1(\mu,Y_t) + \chi_t,  
\enq
for some $\F^0$-adapted processes $(K,\Lambda,Y,\chi)$, with values  in $\S^d\times\S^d\times\R^d\times\R$, and in the backward SDE form
\begin{equation} \label{BSDELQ}
\left\{
\begin{array}{ccll}
dK_t &=& \dot K_t  dt + Z_t^K dW_t^0, \;\;\;  0 \leq t \leq T, &  K_T = P   \\
d\Lambda_t &=& \dot\Lambda_t dt + Z_t^\Lambda dW_t^0, \;\;\; 0 \leq t \leq T, & \Lambda_T = P + \bar P  \\
dY_t &=&  \dot Y_t dt + Z_t^Y dW_t^0, \;\;\; 0 \leq t \leq T, & Y_T = L \\
d\chi_t & =& \dot\chi_t + Z_t^\chi dW_t^0, \;\;\; 0 \leq t \leq T, & \chi_T = 0,
\end{array}
\right.
\end{equation}
for some $\F^0$-adapted processes $\dot K$, $\dot\Lambda$,  $Z^K$, $Z^\Lambda$ with values  in $\S^d$,  $\dot Y$, $Z^Y$ with values in $\R^d$, and $\dot\chi$,  $Z^\chi$ with values in $\R$. Notice that the terminal conditions in \reff{BSDELQ} ensure by \reff{hatfgquadra} that $w$ in  \reff{wLQ} satisfies: $w_T(\mu)$ $=$ $\hat g(\mu)$, and we shall next determine 
the generators $\dot K$, $\dot\Lambda$, $\dot Y$, and $\dot\chi$ in order to satisfy the local (sub)martingale conditions of Lemma \ref{lemverif}. 
Notice that  the functions ${\rm Var}$, $v_2$, $v_1$ are smooth w.r.t. both their arguments, and we have
\begin{equation}\label{derivVar}
\begin{array}{c}
\partial_\mu {\rm Var}(\mu,k)(x)  \;=\;  2 k (x-\bar\mu), \;\; \partial_x\partial_\mu {\rm Var}(\mu,k)(x)  \;=\;   2k \; = \; - \partial_\mu^2 {\rm Var}(\mu,k)(x,x'),  
\\
\partial_k {\rm Var}(\mu,k) \;=\;  {\rm Var}(\mu) \; := \; \int (x-\bar\mu)(x-\bar\mu)\trans \mu(dx)  \\
\partial_\mu v_2(\mu,\ell)(x)  \;=\;  2 \ell \bar\mu, \;\; \partial_x\partial_\mu v_2(\mu,\ell)(x)  \;=\;   0, \;\; \partial_\mu^2 v_2(\mu,\ell)(x,x')  \;=\;  2\ell, \\
\partial_\ell v_2(\mu,\ell) \;=\; \bar\mu\bar\mu\trans \\
\partial_\mu v_1(\mu,y) \; = \; y, \;\;   \partial_x\partial_\mu v_1(\mu,y) \; = \; 0 \; = \;  \partial_\mu^2 v_1(\mu,y)(x,x'), \;\; \partial_y v_1(\mu,y) \; = \; \bar\mu. 
\end{array}
\end{equation}
Let us denote, for any $\alpha$ $\in$ $\Ac$, by $S^\alpha$ the $\F^0$-adapted process equal to $S_t^\alpha$ $=$ 
$w_t(\rho_t^\alpha) + \int_0^t \hat f_s(\rho_s^\alpha,\alpha_s) ds$, $0\leq t\leq T$, and observe then  by It\^o's formula \reff{Ito}  that it is of the form
\beqs
dS_t^\alpha &=& D_t^\alpha dt + \Sigma_t^\alpha dW_t^0, 
\enqs
with a  drift term $D_t^\alpha$ $=$ $\Dc_t(\rho_t^\alpha,\alpha_t,K_t,\Lambda_t,Y_t)$ given by 
\beqs
\Dc_t(\mu,a,k,\ell,y) &=&  \hat f_t(\mu,a)  +  \mu \big( \L_t^{a}  {\rm Var}(\mu,k) +  \L_t^{a}  v_2(\mu,\ell) 
+ \L_t^{a}  v_1(\mu,y)  \big)  \\
& & \;\;\; + \;  \mu\otimes\mu \big( \M_t^{a} {\rm Var}(\mu,k) + 
\M_t^{a} v_2(\mu,\ell)  +  \M_t^{a}  v_1(\mu,y)    \big) \\
& & \; + \; {\rm tr}(\partial_k {\rm Var}(\mu,k)\trans\dot K_t)  + {\rm tr}(\partial_\ell v_2(\mu,\ell)\trans\dot\Lambda_t)   
+ \partial_y v_1(\mu,y)\trans\dot Y_t + \dot\chi_t  \\
& &  + \; {\rm tr}\big(\partial_k \mu \big(\D_t^a {\rm Var}(\mu,k) \big)\trans Z_t^K\big) + 
{\rm tr}\big(\partial_\ell \mu \big(\D_t^a v_2(\mu,\ell) \big)\trans Z_t^\Lambda\big) + 
\partial_y \mu \big(\D_t^a v_1(\mu,\ell) \big)\trans Z_t^Y,
\enqs
for all $t$ $\in$ $[0,T]$, $\mu$ $\in$ $\Pc_{_2}(\R^d)$, $k,\ell$ $\in$ $\S^d$, $y$ $\in$ $\R^d$, $a$ $\in$ $A$. 
(The second-order derivatives terms w.r.t. $k$, $\ell$ and $y$ do not appear since the functions $v_2$, ${\rm Var}$ and $v_1$ are linear, respectively,  in $k$, $\ell$ and $y$, respectively). From the derivatives expression of ${\rm Var}$, $v_2$ and $v_1$ in \reff{derivVar}, we then have
\beq
\Dc_t(\mu,a,k,\ell,y) &=& \hat f_t(\mu,a) + \int  b_t(x,\bar\mu,a)\trans [2k(x-\bar\mu) + 2\ell \bar\mu + y]   \mu(dx) \nonumber \\
& & \;+ \; \int [\sigma_t(x,\bar\mu,a)\trans k \sigma_t(x,\bar\mu,a) + \sigma_t^0(x,\bar\mu,a)\trans k \sigma_t^0(x,\bar\mu,a)] \mu(dx) \nonumber \\
& & \; +\; \Big(\int \sigma_t^0(x,\bar\mu,a)\mu(dx)\Big)\trans (\ell-k) \Big(\int\sigma_t^0(x,\bar\mu,a) \mu(dx) \Big) \nonumber \\
& & \; + \; {\rm Var}(\mu,\dot K_t) + v_2(\mu,\dot\Lambda_t) + v_1(\mu,\dot Y_t) + \dot\chi_t \nonumber \\
& & \; + \;  \int \sigma_t^0(x,\bar\mu,a)\trans [2Z_t^K (x-\bar\mu)  + 2Z_t^\Lambda \bar\mu + Z_t^Y] \mu(dx).   \label{Dc}
\enq
We now distinguish between the cases when the control set $A$ is $\R^m$ (LQCMKV1) or $L(\R^d;\R^m)$ (LQCMKV2).

\subsection{Control set $A$ $=$ $\R^m$} \label{secEuclidian}

From the linear form of $b_t$, $\sigma_t$, $\sigma_t^0$ in \reff{XLQ}, and the quadratic form of $\hat f_t$ in \reff{hatfgquadra}, after some 
straightforward calculations, we have:
\beqs
\Dc_t(\mu,a,k,\ell,y) 
&=&  {\rm Var}(\mu, \Phi_t(k,Z_t^K) + \dot K_t) + v_2(\mu,\Psi_t(k,\ell,Z_t^\Lambda) + \dot\Lambda_t) \\
& &  \; + \;  v_1(\mu, \Theta_t(k,\ell,Z_t^\Lambda,y,Z_t^Y) +\dot Y_t) + \Delta_t(k,\ell,y,Z_t^Y) + \dot\chi_t  \\
& & \; + \; a \trans \Gamma_t(k,\ell) a + [2 U_t\trans(k,\ell,Z_t^\Lambda)\bar\mu  + R_t(k,\ell,y,Z_t^Y)]\trans a 
\enqs 
with
\begin{equation} \label{PhiPsi}
\left\{
\begin{array}{ccl}
\Phi_t(k,Z_t^K) &=& Q_t + B_t\trans k + k B_t + D_t\trans k D_t + (D_t^0)\trans k D_t^0 + (D_t^0)\trans Z_t^K + Z_t^K D_t^0 \\
\Psi_t(k,\ell,Z_t^\Lambda) &=& Q_t + \bar Q_t + (D_t+\bar D_t)\trans k (D_t+\bar D_t) +  (D_t^0+ \bar D_t^0)\trans\ell(D_t^0+\bar D_t^0)     \\
& & \; + \; (B_t+\bar B_t)\trans \ell + \ell(B_t+\bar B_t) + (D_t^0+\bar D_t^0)\trans Z_t^\Lambda + Z_t^\Lambda(D_t^0+\bar D_t^0) \\
\Theta_t(k,\ell,Z_t^\Lambda,y,Z_t^Y) &=& M_t +  (B_t+\bar B_t)\trans y + 2\ell b_t^0 + 2(D_t+\bar D_t)\trans k \gamma_t + 2(D_t^0+\bar D_t^0)\trans \ell \gamma_t^0 \\
&& \; + \; (D_t^0+\bar D_t^0)\trans Z_t^Y + 2 Z_t^\Lambda \gamma_t^0 \\
\Delta_t(k,\ell,y,Z_t^Y) &=& y\trans b_t^0 + \gamma_t\trans k \gamma_t  +  (\gamma_t^0)\trans\ell\gamma_t^0 + (Z_t^Y)\trans \gamma_t^0 \\
\Gamma_t(k,\ell) &=& N_t + F_t\trans k F_t  + (F_t^0)\trans\ell F_t^0 \\
U_t(k,\ell,Z_t^\Lambda) &=&  (D_t+\bar D_t)\trans k F_t + (D_t^0+\bar D_t^0)\trans\ell F_t^0 + \ell C_t + Z_t^\Lambda F_t^0 \\
R_t(k,\ell,y,Z_t^Y) &=& 2 F_t \trans k \gamma_t  + 2 (F_t^0)\trans \ell \gamma_t^0  + C_t\trans y + (F_t^0)\trans Z_t^Y.  
\end{array}
\right.
\end{equation}
Then, after square completion under the condition that $\Gamma_t(k,\ell)$ is positive definite in $\S^m$, 
we  have
\beq
\Dc_t(\mu,a,k,\ell,y) &=&   {\rm Var}(\mu, \Phi_t(k,Z_t^K) + \dot K_t) \nonumber  \\
& & \;+\;   v_2(\mu,\Psi_t(k,\ell,Z_t^\Lambda)  -  U_t(k,\ell,Z_t^\Lambda)\Gamma_t^{-1}(k,\ell)U_t\trans(k,\ell,Z_t^\Lambda)   + \dot\Lambda_t) \nonumber \\
& & \; + \; v_1(\mu, \Theta_t(k,\ell,Z_t^\Lambda,y,Z_t^Y) -  U_t(k,\ell,Z_t^\Lambda)\Gamma_t^{-1}(k,\ell)R_t(k,\ell,y,Z_t^Y)  + \dot Y_t) \nonumber \\
& & \; + \;  \Delta_t(k,\ell,y,Z_t^Y)   -   \frac{1}{4} R_t\trans(k,\ell,y,Z_t^Y)\Gamma_t^{-1}(k,\ell)R_t(k,\ell,y,Z_t^Y)    + \dot\chi_t \nonumber \\
& & \; + \; \big(a - \hat a_t(\bar\mu,k,\ell,y)\big)\trans \Gamma_t(k,\ell) \big(a - \hat a_t(\bar\mu,k,\ell,y)\big),  \nonumber 
\enq
where 
\beqs
\hat a_t(\bar\mu,k,\ell,y) &=& - \Gamma_t^{-1}(k,\ell) \big[ U_t\trans(k,\ell,Z_t^\Lambda) \bar\mu 
+ \frac{1}{2}   R_t(k,\ell,y,Z_t^Y)\big]. 
\enqs
Therefore, whenever 
\beqs
\dot K_t  +  \Phi_t(K_t,Z_t^K)   & =& 0, \\
\dot\Lambda_t +  \Psi_t(K_t,\Lambda_t,Z_t^\Lambda)  
-  U_t(K_t,\Lambda_t,Z_t^\Lambda)\Gamma_t^{-1}(K_t,\Lambda_t)U_t\trans(K_t,\Lambda_t,Z_t^\Lambda)   &=& 0, \\
\dot Y_t + \Theta_t(K_t,\Lambda_t,Z_t^\Lambda,Y_s,Z_t^Y) 
-  U_t(K_t,\Lambda_t,Z_t^\Lambda)\Gamma_t^{-1}(K_t,\Lambda_t)R_t(K_t,\Lambda_t,Y_t,Z_t^Y)  &=& 0, \\
\dot\chi_t + \Delta_t(K_t,\Lambda_t,Y_t,Z_t^Y)   
-   \frac{1}{4} R_t\trans(K_t,\Lambda_t,Y_t,Z_t^Y)\Gamma_t^{-1}(K_t,\Lambda_t,Y_t)R_t(K_t,\Lambda_t,Y_t,Z_t^Y)   &=& 0,
\enqs
holds for all $0\leq t\leq T$,  we have 
\beq
D_t^\alpha &=& \Dc_t(\rho_t^\alpha,\alpha_t,K_t,\Lambda_t,Y_t) \label{Dt} \\
&=& \big(\alpha_t - \hat a_t(\bar\rho_t^\alpha,K_t,\Lambda_t,Y_t)\big)\trans \Gamma_t(K_t,\Lambda_t)
\big(\alpha_t - \hat a_t(\bar\rho_t^\alpha,K_t,\Lambda_t,Y_t)\big), \nonumber 
\enq
which implies that $D_t^\alpha$ $\geq$ $0$, $0\leq t\leq T$, for all $\alpha$ $\in$ $\Ac$, i.e. 
$S_t^\alpha$ $=$ $w_t(\rho_t^\alpha) + \int_0^t \hat f_s(\rho_s^\alpha,\alpha_s) ds$, $0\leq t\leq T$ satisfies the 
$(\P,\F^0)$-local submartingale property for all $\alpha$ $\in$ $\Ac$. We are then led to  consider the system of BSDEs:
\begin{equation} \label{BSDELQ2}
\left\{
\begin{array}{ccl}
dK_t &=& - \Phi_t(K_t,Z_t^K) dt +    Z_t^K dW_t^0, \;\;\;  0 \leq t \leq T,   \; K_T = P \\
d\Lambda_t &=&  - [\Psi_t(K_t,\Lambda_t,Z_t^\Lambda)  
- U_t(K_t,\Lambda_t,Z_t^\Lambda)\Gamma_t^{-1}(K_t,\Lambda_t)U_t\trans(K_t,\Lambda_t,Z_t^\Lambda)] dt  \\
& & \hspace{3cm} + \;  Z_t^\Lambda dW_t^0,  \;\;\;  0 \leq t \leq T, \; \Lambda_T = P + \bar P \\
dY_t &=& -\big[ \Theta_t(K_t,\Lambda_t,Z_t^\Lambda,Y_t,Z_t^Y) -  U_t(K_t,\Lambda_t,Z_t^\Lambda)\Gamma_t^{-1}(K_t,\Lambda_t)
R_t(K_t,\Lambda_t,Y_t,Z_t^Y)] dt \\
& & \hspace{3cm} +  \;   Z_t^Y dW_t^0,  \;\;\;  0 \leq t \leq T, \;\; Y_T = L, \\
d\chi_t &=& -\big[  \Delta_t(K_t,\Lambda_t,Y_t,Z_t^Y)   -   \frac{1}{4} R_t\trans(K_t,\Lambda_t,Y_t,Z_t^Y)\Gamma_t^{-1}(K_t,\Lambda_t)R_t(K_t,\Lambda_t,Y_t,Z_t^Y) \big] dt \\
& & \hspace{3cm} + \;  Z_t^\chi dW_t^0,  \;\;\;  0 \leq t \leq T, \;\; \chi_T = 0. 
\end{array}
\right. 
\end{equation}

\begin{Definition}
A solution to the system of BSDE \reff{BSDELQ2}  is a quadruple of pair $(K,Z^K)$, $(\Lambda,Z^\Lambda)$, $(Y,Z^Y)$, $(\chi,Z^\chi)$ of $\F^0$-adapted processes,  with values, respectively,  in 
$\S^d\times\S^d$, $\S^d\times\S^d$, $\R^d\times\R^d$, $\R\times\R$, respectively, such that 
$\int_0^T |Z_t^K|^2 + |Z_t^\Lambda|^2 + |Z_t^Y|^2 + |Z_t^\chi|^2 dt$ $<$ $\infty$ a.s., 
the matrix process $\Gamma(K,\Lambda)$ with values in  $\S^m$ is positive definite a.s., 
and the following relation 
\begin{equation*} 
\left\{
\begin{array}{ccl}
K_t &=& P +    \int_t^T  \Phi_s(K_s,Z_s^K) ds - \int_t^T  Z_s^K dW_s^0,    \\
\Lambda_t &=& P + \bar P + \int_t^T   \Psi_s(K_s,\Lambda_s,Z_s^\Lambda)  
+ U_s(K_s,\Lambda_s,Z_s^\Lambda)\Gamma_s^{-1}(K_s,\Lambda_s)U_s\trans(K_s,\Lambda_s,Z_s^\Lambda) ds  \\
& & \hspace{3cm} - \; \int_t^T Z_s^\Lambda dW_s^0,   \\
Y_t &=& L+  \int_t^T  \Theta_s(K_s,\Lambda_s,Z_s^\Lambda,Y_s,Z_s^Y) -  U_s(K_s,\Lambda_s,Z_s^\Lambda)\Gamma_s^{-1}(K_s,\Lambda_s)
R_s(K_s,\Lambda_s,Y_s,Z_s^Y) ds \\
& & \hspace{3cm} - \; \int_t^T Z_s^Y dW_s^0,  \\
\chi_t &=& \int_t^T   \Delta_s(K_s,\Lambda_s,Y_s,Z_s^Y)   -   \frac{1}{4} R_s\trans(K_s,\Lambda_s,Y_s,Z_s^Y)\Gamma_s^{-1}(K_s,\Lambda_s)R_s(K_s,\Lambda_s,Y_s,Z_s^Y) ds \\
& & \hspace{3cm} - \; \int_t^T Z_s^\chi dW_s^0, 
\end{array}
\right. 
\end{equation*}
is satisfied for all $t$ $\in$ $[0,T]$.  
\end{Definition}

The following verification result makes the connection between the system \reff{BSDELQ2} and the LQCMKV1 control problem. 

\begin{Proposition} \label{theoverif}
Assume that $(K,Z^K)$, $(\Lambda,Z^\Lambda)$, $(Y,Z^Y)$, $(\chi,Z^\chi)$ is a solution to  BSDE \reff{BSDELQ2} such that 
$K,\Lambda$, $\Gamma^{-1}(K,\Lambda)$ are essentially bounded,  $Z^\Lambda$ lies in $L^2(\Omega\times [0,T])$, i.e., $\E[\int_0^T |Z_t^\Lambda|^2 dt]$ $<$ $\infty$, 
$Y$  lies in $\Sc^2(\Omega\times [0,T])$, i.e. $\E[|\sup_{0\leq t\leq T}|Y_t|^2]$ $<$ $\infty$, and 
$\chi$ lies in $\Sc^1(\Omega\times [0,T])$, i.e.  $\E[|\sup_{0\leq t\leq T}|\chi_t|]$ $<$ $\infty$
Then, the control process
\beq
\alpha_t^* &=& \hat a_t(\E[X_t^*|W^0],K_t,\Lambda_t,Y_t) \label{opticontrol}  \\
&=& - \Gamma_t^{-1}(K_t,\Lambda_t) \big[ U_t\trans(K_t,\Lambda_t,Z_t^\Lambda) \E[X_t^*|W^0]   + \frac{1}{2}   R_t(K_t,\Lambda_t,Y_t,Z_t^Y)\big], \;\;  0 \leq t\leq T,  \nonumber 
\enq
where $X^*$ $=$ $X^{\alpha^*}$   is the state process with the feedback control $\hat a_t(.,K_t,\Lambda_t,Y_t)$, is an optimal control for  the LQCMKV1 problem, i.e., $V_0$ $=$ $J(\alpha^*)$, and we have
\beqs
V_0  &=& {\rm Var}(\Lc(\xi_0),K_0) +  v_2(\Lc(\xi_0),\Lambda_0) + v_1(\Lc(\xi_0),Y_0) + \chi_0. 
\enqs
\end{Proposition}
{\bf Proof.} Consider  $(K,Z^K)$, $(\Lambda,Z^\Lambda)$, $(Y,Z^Y)$, $(\chi,Z^\chi)$  a solution to the BSDE \reff{BSDELQ2}, and $w$ as of the quadratic form \reff{wLQ}. First, notice that $w$ satisfies the quadratic growth \reff{wquadragrowth} since $K,\Lambda$ are essentially bounded, and $(Y,\chi)$ $\in$ $\Sc^2(\Omega\times [0,T])\times\Sc^1(\Omega\times [0,T])$.  Moreover,  we have the terminal condition $w_T(\mu)$ $=$ $\hat g$. 
Next, by construction,  the process $D_t^\alpha$ $=$ $\Dc_t(\rho_t^\alpha,\alpha_t,K_t,\Lambda_t,Y_t)$, $0\leq t\leq T$, 
is nonnegative, which means that $S_t^\alpha$ $=$ $w_t(\rho_t^\alpha) + \int_0^t \hat f_s(\rho_s^\alpha,\alpha_s) ds$, $0\leq t\leq T$, is a $(\P,\F^0)$-local submartingale. Moreover, by choosing the control $\alpha^*$ in the  form \reff{opticontrol}, we notice that $X^*$,  the solution to a linear stochastic McKean-Vlasov dynamics, satisfies the square integrability condition: $\E[\sup_{0\leq t\leq T}|X_t^*|^2]$ $<$ $\infty$, thus 
$\E[\int_0^T |\alpha_t^*|^2 dt]$ $<$ $\infty$,  since $U(K,\Lambda,Z^\Lambda)$ inherits from $Z^\Lambda$ the square integrability condition $L^2(\Omega\times [0,T])$, $\Gamma^{-1}(K,\Lambda)$ is essentially bounded, and so $\alpha^*$  $\in$ 
$\Ac$.  Finally, from \reff{Dt} we see that 
$D^{\alpha^*}$ $=$ $0$, which gives the $(\P,\F^0)$-local martingale property of $S^{\alpha^*}$,  
and we conclude by the dynamic programming verification Lemma \ref{lemverif}. 
\ep

\vspace{2mm}

Let us now show, under assumptions {\bf (H1)} and {\bf (H2)},  the existence of a solution to the BSDE \reff{BSDELQ2} satisfying the integrability  conditions of Proposition \ref{theoverif}. We point out that this system is decoupled: 
\begin{itemize}
\item [(i)] One first considers  
the BSDE for $(K,Z^K)$ whose generator $(k,z)$ $\in$ $\S^d\times\S^d$ $\mapsto$ $\Phi_t(k,z)$ $\in$ $\S^d$ is linear, with essentially bounded coefficients. Since the terminal condition $P$ is also essentially bounded, 
it is known by standard results for linear BSDEs that there exists a unique solution $(K,Z^K)$ with values in $\S^d\times\S^d$, s.t. $K$ is essentially bounded and  $Z^K$ lies in $L^2(\Omega\times [0,T])$.  Moreover, since 
$P$ and $\Phi_t(0,0)$ $=$ $Q_t$ are nonnegative under {\bf (H1)},  we also obtain by standard comparison principle for  BSDE that $K_t$ is nonnegative, for all $0\leq t\leq T$.  
\item[(ii)] Given $K$, we next consider the BSDE for $(\Lambda,Z^\Lambda)$ with generator: $(\ell,z)$ $\in$ $\S^d\times\S^d$ $\mapsto$ 
$\Psi_t(K_t,\ell,z)$ $-$ $U_t(K_t,\ell,z)\Gamma_t^{-1}(K_t,\ell)U_t\trans(K_t,\ell,z)$  $\in$  $\S^d$, and terminal condition $P+\bar P$.  This is a backward stochastic Riccati  equation (BSRE), and it is well-known (see, e.g.,  
\cite{bis76}) that it is  associa\-ted with  a stochastic standard LQ control problem (without McKean-Vlasov dependence) with controlled linear dynamics:  
\beqs
d\tilde X_t &=& [(B_t + \bar B_t)\tilde X_t +  C_t\alpha_t  ] dt + [(D_t^0+\bar D_t^0)\tilde X_t + F_t^0 \alpha_t ] dW_t^0,
\enqs
and quadratic cost functional 
\beqs
\tilde J^K(\alpha) &=& \E \Big[ \int_0^T \big( \tilde X_t\trans Q_t^K \tilde X_t + \alpha_t\trans N_t^K\alpha_t + 2\tilde X_t\trans M_t^K\alpha_t \big) dt \; + \;  \tilde X_T\trans(P +\bar P)\tilde X_T \Big], 
\enqs
where $Q_t^K$ $=$ $Q_t+\bar Q_t + (D_t+\bar D_t)\trans K_t(D_t+\bar D_t)$, $N_t^K$ $=$ $N_t+F_t\trans K_t F_t$, $M_t^K$ $=$ $(D_t+\bar D_t)\trans K_t F_t$.  
Under the condition that $N^K$ is positive definite, we can rewrite this cost functional after square completion as 
\beqs
\tilde J^K(\alpha) &=& \E \Big[ \int_0^T \big( \tilde X_t\tilde Q_t^K \tilde X_t + \tilde\alpha_t\trans N_t^K\tilde\alpha_t \big) dt \; 
+ \;  \tilde X_T\trans(P +\bar P)\tilde X_T \Big],
\enqs
with $\tilde Q_t^K$ $=$ $Q_t^K-M_t^K (N_t^K)^{-1}(M_t^K)\trans$, $\tilde\alpha_t$ $=$  $\alpha_t + (N_t^K)^{-1}(M_t^K)\trans\tilde X_t$. By noting that  $\tilde Q_t^K$ $\geq$ $Q_t+\bar Q_t$, it follows that the symmetric matrices 
$\tilde Q^K$ and $P+\bar P$ are nonnegative under condition {\bf (H1)}, and assuming furthermore that $N^K$ is uniformly positive definite, 
we  obtain  from \cite{tan03} the existence  and uniqueness of a solution $(\Lambda,Z^\Lambda)$ 
to this BSRE, with $\Lambda$ being nonnegative and essentially bounded, and $Z^\Lambda$ square integrable in $L^2(\Omega\times [0,T])$. This implies, in particular,  that $\Gamma^{-1}(K,\Lambda)$ is well-defined and essentially bounded. 
Since $K$ is nonnegative under {\bf (H1)}, notice that the uniform positivity condition on $N^K$ is  satisfied under {\bf (H2)}: this is clear when $N$ is uniformly positive definite (as usually assumed in LQ problem), and holds also true   
when  $F$ is  uniformly nondegenerate,  and $K$ is uniformly positive definite,  which occurs when $P$ or $Q$ is uniformly positive definite  from comparison principle for the linear BSDE for $K$. 
\item[(iii)]  Given $(K,\Lambda,Z^\Lambda)$, we consider the  BSDE for $(Y,Z^Y)$ with generator: $(y,z)$ $\in$ $\R^d\times\R^d$ $\mapsto$ $G_t(y,z)$ $:=$ $\Theta_t(K_t,\Lambda_t,Z_t^\Lambda,y,z)$ $-$ 
$U_t(K_t,\Lambda_t,Z_t^\Lambda)\Gamma_t^{-1}(K_t,\Lambda_t)R_t\trans(K_t,\Lambda_t,y,z)$ with values in $\R^d$, and terminal condition $L$.  This is a linear BSDE and $\{G_t(0,0),0\leq t\leq T\}$ lies in 
$L^2(\Omega\times [0,T])$ (recall that $b^0$, $\gamma$, and $\gamma^0$ are assumed square integrable).  By standard results for BSDEs, we then know that there exists a unique solution $(Y,Z^Y)$ s.t. 
$Y$ lies in $\Sc^2(\Omega\times[0,T])$, and $Z$ lies in $L^2(\Omega\times [0,T])$.  
\item[(iv)] Finally, given $(K,\Lambda,Y,Z^Y)$, we solve the backward stochastic equation for $\chi$, which is explicitly written as 
\beqs
\chi_t &=& \E \Big[   \int_t^T   \Delta_s(K_s,\Lambda_s,Y_s,Z_s^Y)   \\
& & \hspace{1cm} -   \frac{1}{4} R_s\trans(K_s,\Lambda_s,Y_s,Z_s^Y)\Gamma_s^{-1}(K_s,\Lambda_s)R_s(K_s,\Lambda_s,Y_s,Z_s^Y) ds \big| \Fc_t^0 \Big], \; 0 \leq t\leq T,   
\enqs
and $\chi$ satisfies the $\Sc^1(\Omega\times [0,T])$ integrability condition. 
\end{itemize}

To sum up, we have proved  the following result: 

\begin{Theorem} \label{theomain1}
Under assumptions {\bf (H1)} and {\bf (H2)},  
there exists a unique solution $(K,Z^K)$, $(\Lambda,Z^\Lambda)$, $(Y,Z^Y)$, $(\chi,Z^\chi)$  to the BSDE \reff{BSDELQ2} satisfying the 
integrability condition of Proposition \ref{theoverif}, and consequently we have an optimal control for the LQCMKV1 problem given by \reff{opticontrol}.   
\end{Theorem}

\subsection{Control set $A$ $=$ $L(\R^d;\R^m)$}

From the linear form of $b_t$, $\sigma_t$, $\sigma_t^0$ in \reff{XLQ}, and the quadratic form of $\hat f_t$ in \reff{hatfgquadra}, the random field process in \reff{Dc} is given, after some 
straightforward calculations by
\beqs
\Dc_t(\mu,a,k,\ell,y)  &=&  {\rm Var}(\mu, \Phi_t(k,Z_t^K) + \dot K_t) + v_2(\mu,\Psi_t(k,\ell,Z_t^\Lambda) + \dot\Lambda_t) \\
& &  \; + \;  v_1(\mu, \Theta_t(k,\ell,Z_t^\Lambda,y,Z_t^Y) +\dot Y_t) + \Delta_t(k,\ell,y,Z_t^Y) + \dot\chi_t  \\
& & \; + \;  {\rm Var}(a\star\mu,\Gamma_t(k,k)) +  \overline{a\star\mu}\trans \Gamma_t(k,\ell) \overline{a\star\mu}  \\
& & \; + \; 2 \int (x-\bar\mu)\trans  V_t(k,Z_t^K) a(x) \mu(dx) \\
& & \; +  \;  [2 U_t\trans(k,\ell,Z_t^\Lambda)\bar\mu  + R_t(k,\ell,y,Z_t^Y)]\trans \overline{a\star\mu},  
\enqs 
for all $t$ $\in$ $[0,T]$, $\mu$ $\in$ $\Pc_{_2}(\R^d)$, $k,\ell$ $\in$ $\S^d$, $y$ $\in$ $\R^d$, $a$ $\in$ $L(\R^d;\R^m)$,  where $a\star\mu$ $\in$ $\Pc_{_2}(\R^m)$ denotes the image by $a$ of $\mu$,  
\beqs
\overline{a\star\mu} \; = \;  \int a(x) \mu(dx), & &  {\rm Var}(a\star\mu,k) \; = \; \int \Big(a(x) - \overline{a\star\mu} \Big)\trans k \Big(a(x) - \overline{a\star\mu}  \Big) \mu(dx), 
\enqs
and we keep the same notations as in \reff{PhiPsi}  with the additional term:
\beq \label{defV}
V_t(k,Z_t^K) &=& D_t\trans kF_t + (D_t^0)\trans k F_t^0 + kC_t  + Z_t^K F_t^0. 
\enq
Then, after square completion under the condition that $\Gamma_t(k,\ell)$ is positive definite in $\S^m$, we have 
\beqs
\Dc_t(\mu,a,k,\ell,y)  &=&  {\rm Var}(\mu, \Phi_t(k,Z_t^K)  - V_t(k,Z_t^K)\Gamma_t^{-1}(k,k)V_t\trans(k,Z_t^K)   + \dot K_t) \\
& & \;+\;   v_2(\mu,\Psi_t(k,\ell,Z_t^\Lambda)  -  U_t(k,\ell,Z_t^\Lambda)\Gamma_t^{-1}(k,\ell)U_t\trans(k,\ell,Z_t^\Lambda)   + \dot\Lambda_t) \nonumber \\
& & \; + \; v_1(\mu, \Theta_t(k,\ell,Z_t^\Lambda,y,Z_t^Y) -  U_t(k,\ell,Z_t^\Lambda)\Gamma_t^{-1}(k,\ell)R_t(k,\ell,y,Z_t^Y)  + \dot Y_t) \nonumber \\
& & \; + \;  \Delta_t(k,\ell,y,Z_t^Y)   -   \frac{1}{4} R_t\trans(k,\ell,y,Z_t^Y)\Gamma_t^{-1}(k,\ell)R_t(k,\ell,y,Z_t^Y)    + \dot\chi_t \nonumber  \\
& & \; + \; {\rm Var}\big( (a-\hat\ba_t)(.,\bar\mu,k,\ell,y)\star\mu,\Gamma_t(k,k)\big) \\
& & \; + \;  \overline{(a-\hat\ba_t)(.,\bar\mu,k,\ell,y)\star\mu}\trans \Gamma_t(k,\ell) \overline{(a-\hat\ba_t)(.,\bar\mu,k,\ell,y)\star\mu}
\enqs
where $\hat \ba_t(.,\bar\mu,k,\ell,y)$ $:$ $\R^d$ $\rightarrow$ $\R^m$ is defined by
\beqs
\hat \ba_t(x,\bar\mu,k,\ell,y) &=& - \Gamma_t^{-1}(k,k)V_t(k,Z_t^K)\trans (x-\bar\mu) \\
& & \;\;\;\;\;  - \;   \Gamma_t^{-1}(k,\ell) \big[ U_t\trans(k,\ell,Z_t^\Lambda) \bar\mu  + \frac{1}{2}   R_t(k,\ell,y,Z_t^Y)\big], \;\;\; x \in \R^d. 
\enqs
We then consider the system of BSDEs:
\begin{equation} \label{BSDELQ3}
\left\{
\begin{array}{ccl}
dK_t &=& - \big[\Phi_t(K_t,Z_t^K)  -  V_t(K_t,Z_t^K)\Gamma_t^{-1}(K_t,K_t)V_t\trans(K_t,Z_t^K) \big]dt \\
& &\hspace{3cm}  + \;    Z_t^K dW_t^0, \;\;\;  0 \leq t \leq T,   \; K_T = P \\
d\Lambda_t &=&  - \big[\Psi_t(K_t,\Lambda_t,Z_t^\Lambda)  
- U_t(K_t,\Lambda_t,Z_t^\Lambda)\Gamma_t^{-1}(K_t,\Lambda_t)U_t\trans(K_t,\Lambda_t,Z_t^\Lambda)\big] dt  \\
& & \hspace{3cm} + \;  Z_t^\Lambda dW_t^0,  \;\;\;  0 \leq t \leq T, \; \Lambda_T = P + \bar P \\
dY_t &=& -\big[ \Theta_t(K_t,\Lambda_t,Z_t^\Lambda,Y_t,Z_t^Y) -  U_t(K_t,\Lambda_t,Z_t^\Lambda)\Gamma_t^{-1}(K_t,\Lambda_t)
R_t(K_t,\Lambda_t,Y_t,Z_t^Y)] dt \\
& & \hspace{3cm} +  \;   Z_t^Y dW_t^0,  \;\;\;  0 \leq t \leq T, \;\; Y_T = L, \\
d\chi_t &=& -\big[  \Delta_t(K_t,\Lambda_t,Y_t,Z_t^Y)   -   \frac{1}{4} R_t\trans(K_t,\Lambda_t,Y_t,Z_t^Y)\Gamma_t^{-1}(K_t,\Lambda_t)R_t(K_t,\Lambda_t,Y_t,Z_t^Y) \big] dt \\
& & \hspace{3cm} + \;  Z_t^\chi dW_t^0,  \;\;\;  0 \leq t \leq T, \;\; \chi_T = 0, 
\end{array}
\right. 
\end{equation}
and by the same arguments as in Proposition \ref{theoverif}, we have the following verification result making the connection between the system \reff{BSDELQ3} and the LQCMKV2 control problem. 

\begin{Proposition} \label{theoverif2}
Assume that $(K,Z^K)$, $(\Lambda,Z^\Lambda)$, $(Y,Z^Y)$, $(\chi,Z^\chi)$ is a solution to the BSDE \reff{BSDELQ3} such that 
$K,\Lambda$, $\Gamma^{-1}(K,\Lambda)$ are essentially bounded, $Y$  lies in $\Sc^2(\Omega\times [0,T])$,  and $\chi$ lies in $\Sc^1(\Omega\times [0,T])$. 
Then, the control process $\alpha^*$ with values in $L(\R^d;\R^m)$ and defined by  
\beq
\alpha_t^*(x)  &=& \hat \ba_t(x,\E[X_t^*|W^0],K_t,\Lambda_t,Y_t) \label{opticontrol2}  \\
&=&  - \Gamma_t^{-1}(K_t,K_t)V_t(K_t,Z_t^K)\trans (x- \E[X_t^*|W^0])  \nonumber \\
& &  - \;  \Gamma_t^{-1}(K_t,\Lambda_t) \big[ U_t\trans(K_t,\Lambda_t,Z_t^\Lambda) \E[X_t^*|W^0]   + \frac{1}{2}   R_t(K_t,\Lambda_t,Y_t,Z_t^Y)\big], \;  x \in \R^d, \; 0 \leq t\leq T,  \nonumber 
\enq
where $X^*$ $=$ $X^{\alpha^*}$   is the state process with the feedback control $\hat\ba_t(.,.,K_t,\Lambda_t,Y_t)$, is an optimal control for  the LQCMKV2 problem, i.e.,  $V_0$ $=$ $J(\alpha^*)$, and we have
\beqs
V_0  &=& {\rm Var}(\Lc(\xi_0),K_0) +  v_2(\Lc(\xi_0),\Lambda_0) + v_1(\Lc(\xi_0),Y_0) + \chi_0. 
\enqs
\end{Proposition}

\vspace{2mm}

Let us now discuss the existence  of a solution to the BSDE \reff{BSDELQ3} satisfying the integrability  conditions of Proposition \ref{theoverif2}.  As for \reff{BSDELQ2}, this system is decoupled. 
The difference w.r.t to the LQCMKV1 problem is in the BSDE for $(K,Z^K)$,  where the generator  $(k,z)$ $\in$ $\S^d\times\S^d$  $\mapsto$ 
$\Phi_t(k,z)-V_t(k,z)\Gamma_t^{-1}(k,k)V_t\trans(k,z)$ $\in$ $\S^d$ is now of the  Riccati type.  In general, it is not in the class of BSREs related to LQ control problem, but existence can be obtained in some particular cases:    
\begin{itemize}
\item[(1)] The coefficients $B$, $C$, $D$, $F$, $D^0$, $F^0$,  $Q$, $P$,  $N$ are deterministic. In this case, the BSRE for $K$ is reduced to a matrix Riccati ordinary differential equation:
\beqs
- \frac{dK_t}{dt} &=& \Phi_t(K_t,0) - V_t(K_t,0) \Gamma_t^{-1}(K_t,K_t)V_t\trans(K_t,0), \;\;\; 0 \leq t \leq T, \; K_T \; = \; P. 
\enqs
This problem is associated to the LQ problem with controlled linear dynamics 
\beqs
d\tilde X_t &=& (B_t \tilde X_t + C_t \tilde\alpha_t) dt +  (D_t\tilde X_t + F_t\tilde\alpha_t) dW_t + (D_t^0\tilde X_t + F_t^0\tilde\alpha_t) dW_t^0, 
\enqs
where the control process  $\tilde\alpha$ is an $\F$-adapted process  with values in $\R^m$, and the cost functional to be minimized over $\tilde\alpha$ is 
\beqs
\tilde J(\tilde\alpha) &=& \E \Big[ \int_0^T (\tilde X_t\trans Q_t \tilde X_t + \tilde\alpha_t\trans N_t \tilde\alpha_t) dt + \tilde X_T\trans P \tilde X_T \Big].
\enqs
It was solved in \cite{won68} under assumption {\bf (H1)} and the condition {\bf (H2)}(i) that $N$ is uniformly positive definite,  and this gives the existence and uniqueness of $K$ $\in$ $C^1([0,T];\S^d)$, which is nonnegative.  
\item[(2)]  $D$ $\equiv$ $F$ $\equiv$ $0$. In this case, the BSRDE for $(K,Z^K)$ is associated to the LQ problem with controlled linear dynamics 
\beqs
d\tilde X_t &=& (B_t \tilde X_t + C_t \tilde\alpha_t) dt +  (D_t^0\tilde X_t + F_t^0\tilde\alpha_t) dW_t^0, 
\enqs
where the control process  $\tilde\alpha$ is an $\F^0$-adapted process  with values in $\R^m$, and the cost functional to be minimized over $\tilde\alpha$ is 
\beqs
\tilde J(\tilde\alpha) &=& \E \Big[ \int_0^T (\tilde X_t\trans Q_t \tilde X_t + \tilde\alpha_t\trans N_t \tilde\alpha_t) dt + \tilde X_T\trans P \tilde X_T \Big].
\enqs
It is then known from \cite{tan03} that under assumptions {\bf (H1)} and  {\bf (H2)}(i),  there exists a unique pair $(K,Z^K)$ solution to the BSRDE, with $K$ nonnegative, and essentially bounded. 
\item[(3)] $N$ $\equiv$ $0$, $P$ is uniformly positive, $m$ $=$ $d$, and $F$ is invertible with $F^{-1}$ bounded. 
In this case, the BSDE for $K$ is reduced to the linear BSDE:
\beqs
dK_t &=& - \big[ \Phi_t(K_t,Z_t^K) - (C_t F_t^{-1} D_t)\trans K_t + K_t(C_t F_t^{-1} D_t)  \\
& & \;\;  - \;   D_t\trans K_t D_t - K_t C_t(F_t\trans K_t F_t)^{-1} C_t\trans K_t  \big] dt  + Z_t^K dW_t^0, \;\; 0 \leq t\leq T, \; K_T = P,
\enqs
for which it is known that there exists a unique solution  $(K,Z^K)$, with $K$ positive, and essentially bounded. 
\end{itemize}
It is an open question whether existence of a solution for $K$ to the BSRE \reff{BSDELQ3} holds  in the general case. Anyway, once a solution $K$ exists, and is given, the BSDEs for the pairs $(\Lambda,Z^\Lambda)$, $(Y,Z^Y)$, $(\chi,Z^\chi)$ are the same as in \reff{BSDELQ2}, and then  their existence and uniqueness are  obtained under the same conditions.

\section{Applications}

\setcounter{equation}{0} \setcounter{Assumption}{0}
\setcounter{Theorem}{0} \setcounter{Proposition}{0}
\setcounter{Corollary}{0} \setcounter{Lemma}{0}
\setcounter{Definition}{0} \setcounter{Remark}{0}

\subsection{Trading with price impact and benchmark tracking
}



We consider an agent trading in a financial market  with an inventory $X_t$, i.e., a  number of shares held at time $t$ in a risky stock,  governed by 
\beqs
dX_t &=& \alpha_t dt,
\enqs
where the control $\alpha$,  a real-valued $\F^0$-progressively measurable process in  $L^2(\Omega\times [0,T])$, represents the trading rate. 
Given a real-valued $\F^0$-adapted stock price process $(S_t)_{0\leq t\leq T}$ in  $L^2(\Omega\times [0,T])$, a real-valued $\F^0$-adapted target process $(I_t)_{0\leq t\leq T}$ in  $L^2(\Omega\times [0,T])$, and 
a terminal benchmark $H$ as a square integrable  $\Fc_T^0$-measurable random variable, 
the objective of the agent is to minimize over control processes $\alpha$ a cost functional  of the form: 
\beq \label{Jcostgen}
J(\alpha) &=& \E \Big[ \int_0^T \Big( \alpha_t \big( S_t + \eta \alpha_t)  + q(X_t -  I_t)^2  \Big) dt + \lambda (X_T - H)^2 \Big],     
\enq
where $\eta$ $>$ $0$, $q$ $\geq$ $0$, and $\lambda$ $\geq$ $0$ are constants. 

Such  formulation is connected with optimal trading and hedging problems in presence of liquidity frictions like price impact, and widely studied in the recent years: when $S$ $\equiv$ $=$ $0$, 
the cost functional in \reff{Jcostgen} arises in  option hedging in presence of transient price impact, see, e.g.,   
\cite{rogsin10}, \cite{almli15}, \cite{bansonvos15}, and is also related to the problem of optimal VWAP execution (see \cite{frewes13}, \cite{carjai15}), 
or benchmark tracking, see \cite{cairostan16}.  When $q$ $=$ $0$,  the minimization of  the cost functional in \reff{Jcostgen} corresponds to the  optimal execution problem arising in limit order book (LOB), as originally formulated in \cite{almcri00} in a particular Bachelier model for $S$, and has been extended (with general shape functions in LOB) in the literature, but mostly  by assuming the martingale property of the price process, see, e.g.,  
\cite{alfetal}, \cite{preetal}.  By rewriting  the cost functional  after square completion as
\beqs
J(\alpha) &=&   \E \Big[ \int_0^T  \big(  \eta \tilde\alpha_t^2  + q(X_t -  I_t)^2  \big) dt +  \lambda (X_T - H)^2 \Big] \; - \;  \E\Big[ \int_0^T \frac{S_t^2}{4\eta} dt \Big],
\enqs
with $\tilde\alpha_t$ $=$ $\alpha_t$ $+$ $\frac{S_t}{2\eta}$, we see that this problem fits into the LQCMKV1 framework (with $b_t^0$ $=$ 
$- \frac{S_t}{2\eta}$, without McKean-Vlasov dependence but with random coefficients), and Assumptions {\bf (H1)}, {\bf (H2)} are satisfied. 
From Theorem \ref{theomain1}, the optimal control is then given by 
\beq \label{alphainter}
\alpha_t^* &=& -  \frac{1}{\eta} \big[\Lambda_t X_t^* + \frac{Y_t}{2} \big] - \frac{S_t}{2\eta}, \;\;\; 0 \leq t\leq T, 
\enq
where $\Lambda$ is solution to the  (ordinary differential) Riccati equation
\beq \label{ric1}
d\Lambda_t &=& -(q - \frac{\Lambda_t^2}{\eta}) dt, \;\;\; 0 \leq t\leq T, \; \Lambda_T \; = \; \lambda, 
\enq
and $Y$ is solution to the linear BSDE
\beq \label{Y1}
dY_t &=&  
\big[ 2qI_t   +  \frac{\Lambda_t}{\eta} S_t +   \frac{\Lambda_t}{\eta} Y_t  \big] dt +  Z_t^Y dW_t^0, \;\;\; 0 \leq t\leq T, \; Y_T \; = \; -2\lambda H.  
\enq
The solution to the Riccati equation is 
\beqs
\frac{\Lambda_t}{\eta} &=& \sqrt{q/\eta}  \frac{ \sqrt{q/\eta}\sinh(\sqrt{q/\eta}(T-t)) + \lambda/\eta\cosh(\sqrt{q/\eta}(T-t))}{\lambda/\eta \sinh(\sqrt{q/\eta}(T-t)) + \sqrt{q/\eta}\cosh(\sqrt{q/\eta}(T-t))}, \;\;\; 0 \leq t \leq T, 
\enqs 
while the solution to the linear BSDE is given by
\beqs
Y_t &=& 
- 2 \E \Big[ e^{-\int_t^T \frac{\Lambda_s}{\eta} ds} \lambda H + \int_t^T e^{-\int_t^s \frac{\Lambda_u}{\eta} du} \big(q I_s + \frac{\Lambda_s}{\eta} S_s \big)   ds \big| \Fc_t^0 \Big], 
\;\; 0 \leq t \leq T. 
\enqs 
By integrating the function $\Lambda/\eta$, we have 
\beqs
e^{-\int_t^s \frac{\Lambda_u}{\eta} du} &=&  \frac{\Lambda_t/\eta}{\sqrt{q/\eta}} \frac{ \sqrt{q/\eta} \cosh(\sqrt{q/\eta}(T-s)) + \lambda/\eta\sinh(\sqrt{q/\eta}(T-s) )} {  \sqrt{q/\eta} \sinh(\sqrt{q/\eta}(T-t)) + \lambda/\eta\cosh(\sqrt{q/\eta}(T-t))}   \\
&=&   \frac{\Lambda_t}{\Lambda_s} \frac{\sqrt{q/\eta} \sinh(\sqrt{q/\eta}(T-s)) + \lambda/\eta\cosh(\sqrt{q/\eta}(T-s) )} 
{ \sqrt{q/\eta} \sinh(\sqrt{q/\eta}(T-t)) + \lambda/\eta\cosh(\sqrt{q/\eta}(T-t))}, \;\;\; t \leq s \leq T,
\enqs
and plugging into the expectation form of $Y$,  the optimal control in \reff{alphainter}  is then expressed as 
\beq 
\alpha_t^* &=& - \frac{\Lambda_t}{\eta} ( X_t^* -  \hat I_t^H)    + \frac{1}{2\eta} \Big(   \E\big[\int_t^T  \frac{\Lambda_t}{\eta}\frac{\omega(t,T)}{\omega(s,T)} S_s ds | \Fc_t^0\big] - S_t \Big) \nonumber \\
& =: & \alpha_t^{*,IH} + \alpha_t^{*,S}, \;\;\; 0 \leq t \leq T,   \label{alphaoptgen}
\enq
where 
\beqs
\hat I_t^H &=& \E \Big[ \omega(t,T)   H   +  (1- \omega(t,T)) \int_t^T  I_s  \Kc(t,s) ds  
\big| \Fc_t^0 \Big]
\enqs
with a weight  valued in $[0,1]$ 
\beqs
\omega(t,T) &=& \frac{\lambda/\eta}{ \sqrt{q/\eta} \sinh(\sqrt{q/\eta}(T-t)) + \lambda/\eta\cosh(\sqrt{q/\eta}(T-t))},
\enqs
and a kernel
\beqs
\Kc(t,s) &=&  \sqrt{q/\eta} \frac{\sqrt{q/\eta} \cosh(\sqrt{q/\eta}(T-t)) + \lambda/\eta\sinh(\sqrt{q/\eta}(T-t))}{\sqrt{q/\eta} \sinh(\sqrt{q/\eta}(T-t)) + \lambda/\eta(\cosh(\sqrt{q/\eta}(T-t))-1)}, \;\;\;   0\leq t \leq s \leq T. 
\enqs

The optimal trading rule in \reff{alphaoptgen} is decomposed in two parts: 
\begin{itemize}
\item[(i)] The first term $\alpha^{*,IH}$ prescribes the agent to trade 
optimally towards a weigh\-ted ave\-rage $\hat I_t^H$,  rather than the current target position $I$. Indeed, $\hat I^H$ is a convex combination of the expected future of the  terminal random target $H$, and of a weighted average of the running target $I$ (notice that $\Kc(t,.)$ is a nonnegative kernel integrating to one over $[t,T])$.  The rate towards this target is at a speed proportional to its  distance w.r.t the current investor's position, and the coefficient of proportionality is determined by the costs parameters $\eta$, $q$, $\lambda$ and the time to maturity $T-t$. We retrieve the interpretation and results  obtained in  
\cite{bansonvos15} in the limiting cases where $\lambda$ $=$ $0$ (no constraint on the terminal position), and $\lambda$ $=$ $\infty$ (constraint on the terminal position $X_T$ $=$ $H$).  
In the case where $q$ $=$ $0$,  we have $\Lambda_t/\eta$ $=$ $\lambda/(\eta+\lambda(T-t))$,   $\hat I_t^H$ $=$  $\E[H|\Fc_t^0]$, and we retrieve, in particular,  the expression $\alpha^{*,IH}$ $=$ 
$-X_t^*/(T-t)$, of optimal trading rate when $H$ $=$ $0$, and  $\lambda$ $\rightarrow$ $\infty$ corresponding to the optimal execution problem with terminal liquidation $X_T$ $=$ $0$.  
\item[(ii)]  The second term $\alpha^{*,S}$ related to the stock price, is an incentive to buy or sell   depending on whether the weighted average of  expected future value of the stock is larger or smaller than its current value. In particular, 
when  the price process is a martingale,  then 
\beqs
\alpha_t^{*,S} & = & - \frac{S_t}{2\eta} \frac{\sqrt{q/\eta}}{\sqrt{q/\eta} \cosh(\sqrt{q/\eta}(T-t)) + \lambda/\eta\sinh(\sqrt{q/\eta}(T-t))}
\enqs
which is nonpositive for nonnegative price $S_t$, hence meaning that due to the price impact, one must sell. Moreover, in the limiting case where $\lambda$ $\rightarrow$ $\infty$, i.e., the terminal inventory $X_T$ is constrained to achieve the target $H$, then $\alpha^{*,S}$ is zero: we retrieve the result that the optimal trading rate does not depend on the price process when it is a martingale, see \cite{alfetal}, \cite{preetal}. 
\end{itemize} 
On the other hand, by applying It\^o's formula to \reff{alphainter}, and using \reff{ric1}-\reff{Y1}, we have
\beqs
d\big(\alpha_t^* + \frac{S_t}{2\eta} \big)  &=&  \frac{q}{\eta}(X_t^*-I_t)  ds - \frac{1}{2\eta} Z_t^Y dW_s^0, 
\enqs
which implies the notable property: 
\beqs
\alpha_t^* + \frac{S_t}{2\eta}  - \frac{q}{\eta}  \int_0^t (X_s^*-I_s) ds, \;\;\; 0 \leq t\leq T, &  & \mbox{ is a martingale.}  
\enqs


\subsection{Conditional Mean-variance portfolio selection in incomplete market}

We consider  an agent who can invest in a financial market model with one bond of price process $S^0$ and one  risky asset of price process $S$ governed by
\beqs
dS_t^0 &=& S_t^0 r(I_t) dt \\
dS_t &=& S_t ( (b+r)(I_t) dt + \sigma(I_t) dW_t), 
\enqs
 where $I$ is a factor process with dynamics governed by a Brownian motion $W^0$, assumed to be non correlated with the Brownian motion $W$ driving the asset price process $S$, and $r$ the interest rate, $b$ the excess rate of return, and  $\sigma$ the volatility are measurable bounded functions of $I$, with $\sigma(I_t)$ $\geq$ $\eps$ for some $\eps$ $>$ $0$.  
 We shall assume that the  natural filtration generated by the observable  factor process $I$ is equal to the filtration $\F^0$ generated by $W^0$.  Notice that the market is incomplete as the agent cannot trade in the factor process. 
The investment strategy of the agent is modeled by a random field $\F^0$-progressively measurable process $\alpha$ $=$ 
 $\{\alpha_t(x), 0 \leq t\leq T, x  \in \R \}$ (or equivalently as a $\F^0$-progressively measurable process with values in $L(\R;\R)$) where $\alpha_t(x)$ with values in $\R$,  is Lipschitz in $x$, and represents the amount invested in the stock at time $t$, when the current wealth is $X_t$ $=$ $x$, and based on the past observations  $\Fc_t^0$ of the factor process.  The evolution of the controlled wealth process is then given by
 \beq \label{dynXrich}
 dX_t &=&   r(I_t) X_t dt  +  \alpha_t(X_t) \big( b(I_t) dt + \sigma(I_t) dW_t \big), \;\;\;  0 \leq t \leq T, \; X_0 = x_0 \in \R.  
 \enq
 The objective of the agent is to minimize over investment strategies a criterion  of the form:
 \beqs
 J(\alpha) &=& \E \Big[ \frac{\lambda}{2} {\rm Var}(X_T|W^0)  - \E[X_T|W^0] \Big],   
 \enqs
 where $\lambda$ is a  positive $\Fc_T^0$-measurable random variable. In the absence of random factors in the dynamics of the price process,  hence in a complete market model, and when $\lambda$ is constant,  the above criterion reduces to the classical mean-variance portfolio selection, as studied e.g. in \cite{lizho00}.  Here, in presence of the random factor, we consider the expectation of a conditional mean-variance criterion, and also allow  the  risk-aversion parameter 
 $\lambda$ to depend reasonably on the random factor environment.   By rewriting the cost functional as
 \beqs
 J(\alpha) &=& \E \Big[ \frac{\lambda}{2} X_T^2 -  \frac{\lambda}{2} \big( \E[X_T|W^0] \big)^2 - X_T \Big], 
 \enqs
we then see that this conditional mean-variance portfolio selection  problem  fits into the LQCMKV2 problem, and more specifically into the case (3) 
of the discussion following Proposition \ref{theoverif2}. The optimal control is then given from \reff{opticontrol2} by
\beq \label{optialphaMV}
\alpha_t^*(x) &=& - \frac{b(I_t)}{\sigma^2(I_t)}\big(x -  \E[X_t^*|W^0]\big) 
- \frac{b(I_t)}{\sigma^2(I_t) K_t} \big[ \Lambda_t \E[X_t^*|W^0]  + \frac{1}{2} Y_t \big],
\enq
where $X^*$ is the optimal wealth process in \reff{dynXrich} controlled by $\alpha^*$,  $K$ is the solution to the linear BSDE
\beqs
dK_t &=& \big[ \frac{b^2(I_t)}{\sigma^2(I_t)} - 2 r(I_t) \big]K_t  dt  + Z_t^K dW_t^0, \;\;\; 0 \leq t \leq T, \; K_T = \frac{\lambda}{2},
\enqs
$\Lambda$ is solution to the linear BSDE
\beqs
d\Lambda_t &=& \big[  \frac{b^2(I_t)}{\sigma^2(I_t)K_t} \Lambda_t^2 -   2 r(I_t) \Lambda_t \big] dt 
+ Z_t^\Lambda dW_t^0, \;\;\; 0 \leq t\leq T, \; \Lambda_T = 0, 
\enqs
and $Y$ the solution to the linear BSDE
\beqs
dY_t &=& \big[ \frac{b^2(I_t)\Lambda_t}{\sigma^2(I_t)K_t} - r(I_t) \big] Y_t dt + Z_t^Y dW_t^0, \;\;\; 0 \leq t \leq T, \; Y_T = -1. 
\enqs
The solutions to these linear BSDEs are explicitly given by 
\beq \label{Kexpli}
K_t &=& \E \Big[ \frac{\lambda}{2} \exp\Big( \int_t^T   2 r(I_s) - \frac{b^2(I_s)}{\sigma^2(I_s)}   ds\Big) \big| \Fc_t^0 \Big],   
\enq
$\Lambda$ $=$ $0$, and 
\beq \label{Yexpli}
Y_t &=& - \E \Big[ \exp\big( \int_t^T r(I_s) ds \big) \big| \Fc_t^0 \Big], \;\;\; 0 \leq t \leq T. 
\enq 
From \reff{dynXrich} and \reff{optialphaMV}, the conditional mean of the optimal wealth process $X^*$ with portfolio strategy $\alpha^*$ is governed by
\beqs
d \E[X_t^*|W^0] &=& \big[  r(I_t)  \E[X_t^*|W^0]   -  \frac{b^2(I_t)}{2\sigma^2(I_t)} \frac{Y_t}{K_t} \big] dt, 
\enqs
hence explicitly given by
\beqs
\E[X_t^*|W^0] &=&  x_0 e^{\int_0^t r(I_s) ds} 
- \int_0^t  \frac{b^2(I_s)}{2\sigma^2(I_s)} \frac{Y_s}{K_s} e^{\int_s^t r(I_u) du} ds, \;\;\; 0 \leq t \leq T. 
\enqs
Plugging into \reff{optialphaMV}, this gives the explicit form of the optimal control for the conditional mean-variance portfolio selection problem:
\beq \label{optialphagen}
\alpha_t^*(X_t^*) &=& \frac{b(I_t)}{\sigma^2(I_t)} \Big[ x_0  e^{\int_0^t r(I_s) ds}  - X_t^* 
+ \frac{1}{2}\big(\int_0^t   \frac{b^2(I_s)}{\sigma^2(I_s)} \frac{|Y_s|}{K_s} e^{\int_s^t r(I_u) du} ds +  \frac{|Y_t|}{K_t} \big) \Big],  
\enq
for all $0\leq t\leq T$, with $K$ and $Y$ in \reff{Kexpli}-\reff{Yexpli}.  When $b$, $\sigma$, and $r$ do not depend on $I$, we retrieve the expression of the optimal control obtained in \cite{lizho00}, and the formula \reff{optialphagen} is an extension to the case of an incomplete market with a factor $I$ independent of the stock price.

\subsection{Systemic risk model}

We consider a model of inter-bank borrowing and lending  where the log-monetary reserves $X^i$, $i$ $=$ $1,\ldots,n$, of $n$ banks are  driven by
\beqs
dX_t^i &=& \frac{\kappa(I_t)}{n} \sum_{j=1}^n (X_t^j - X_t^i) dt  + \alpha_t^i dt + \sigma(I_t) (\sqrt{1-\rho^2(I_t)} dW_t^i + \rho(I_t) dW_t^0),  
\; i=1,\ldots,n, 
\enqs
where $I_t$ is a factor process  driven by  a Brownian motion $W^0$, which is the common noise for all the banks, $W^i$, $i$ $=$ $1,\ldots,N$, are independent Brownian motions,  independent of $W^0$, called idiosyncratic noises, 
$\rho(I_t)$ $\in$ $[-1,1]$ is the correlation between the idiosyncratic noise and the common noise, $\kappa(I_t)$ $\geq$ $0$ is  the rate of mean-reversion in the interaction from borrowing and lending between the banks,  $\sigma(I_t)$ $>$ $0$ is the volatility of the bank reserves, and compared to the original model introduced in  \cite{caretal14},  these coefficients may depend on the common factor process $I$. Each bank $i$ can control its rate  of borrowing/lending to a central bank via the control $\alpha_t^i$ in order to minimize 
\beqs
J^i(\alpha^1,\ldots,\alpha^n) &=& \E \Big[ \int_0^T  f_t(X_t^i, \frac{1}{n} \sum_{j=1}^n X_t^j,\alpha_t^i) dt + g(X_T^i,\frac{1}{n} \sum_{j=1}^n X_t^j) \Big],
\enqs  
where 
\beqs
f_t(x,\bar x,a) \; = \;  \frac{1}{2} a^2 -  q(I_t) a(x-\bar x) + \frac{\eta(I_t)}{2}(x-\bar x)^2, & & g(x,\bar x) \; = \;  \frac{c}{2} (x-\bar x)^2.
\enqs
Here  $q(I_t)$ $>$ $0$ is a positive $\F^0$-adapted process for the incentive to borrowing ($\alpha_t^i$ $>$ $0$) or lending ($\alpha_t^i$ $<$ $0$),  
$\eta(I_t)$ $>$ $0$ is a positive $\F^0$-adapted process,  $c$ $>$ $0$ is a  positive $\Fc_T^0$-measurable random variable,  for penalizing departure from the average, and these coefficients  may depend on the random factor.   
For this $n$-player stochastic differential game, one looks for cooperative equilibriums by taking the point of view of  
a center of decision (or social planner), which decides on the strategies for all banks, with the goal of minimizing the global cost to the collective.   
More precisely, given the symmetry of the set-up, when the social planner  chooses  the same control policy for all the banks in feedback form: $\alpha_t^i$ $=$ $\tilde\alpha(t,X_t^i,\frac{1}{n}\sum_{j=1}^n X_t^j,I_t)$, $i$ $=$ $1,\ldots,n$, for some deterministic function  $\tilde\alpha$  depending upon time, private state of bank $i$,  the empirical mean of all banks, and factor $I$,  then the theory of propagation of chaos implies that, in the limit $n$ $\rightarrow$ $\infty$, the log-monetary reserve processes $X^i$ become asymptotically independent conditionally on the random environment $W^0$, and the empi\-rical mean $\frac{1}{n}\sum_{j=1}^n X_t^j$ converges to the conditional mean $\E[X_t|W^0]$ of $X_t$ given $W^0$,  and  $X$ is governed by the conditional McKean-Vlasov equation: 
\beqs
dX_t &=&  \big[ \kappa(I_t)(\E[X_t | W^0] - X_t)  + \tilde\alpha(t,X_t,\E[X_t|W^0],I_t) ] dt  \nonumber \\
& & \;\;\; + \;  \sigma(I_t)(\sqrt{1-\rho^2(I_t)} dW_t + \rho(I_t) dW_t^0), \; X_0 \; = \;  x_0  \in   \R,     \label{logX}
\enqs 
for some Brownian motion $W$ independent of $W^0$. More generally, the representative bank can control its rate  of borrowing/lending via a random field $\F^0$-adapted process $\alpha$ $=$ $\{\alpha_t(x),x\in \R\}$, leading to the  log-monetary reserve dynamics:
\beq
dX_t &=&  \big[ \kappa(I_t)(\E[X_t | W^0] - X_t)  + \alpha_t(X_t) ] dt  \nonumber \\
& & \;\;\; + \;  \sigma(I_t)(\sqrt{1-\rho^2(I_t)} dB_t + \rho(I_t) dW_t^0), \; X_0 \; = \;  x_0  \in   \R,     \label{logX1}
\enq 
and the  objective is to minimize over $\alpha$
\beqs
J(\alpha) &=& \E \Big[ \int_0^T  f_t(X_t,\E[X_t|W^0],\alpha_t(X_t)) dt + g(X_T,\E[X_T|W^0]) \Big].
\enqs
After square completion, we can rewrite the cost functional as
\beqs
J(\alpha) &=&  \E \Big[ \int_0^T \Big( \frac{1}{2} \bar\alpha_t(X_t)^2   + \frac{(\eta-q^2)(I_t)}{2} (\E[X_t| W^0] - X_t)^2 \Big) dt    +    \frac{c}{2}  (\E[X_T | W^0] - X_T)^2 \Big],   
\enqs
with $\bar\alpha_t(X_t)$ $=$ $\alpha_t(X_t) - q(\E[X_t|W^0]-X_t)$.  Assuming that $q^2$ $\leq$ $\eta$, 
this model fits into the LQCMKV2  problem, and more specifically into the case (2)  
of the discussion following Proposition \ref{theoverif2}. The optimal control is then given from \reff{opticontrol2} by
\beq \label{optiinter3}
\alpha_t^*(x) &=& -(2K_t +q(I_t))  (x - \E[X_t^*|W^0])  - 2 \Lambda_t \E[X_t^*|W^0] - Y_t, \;  
\enq
where $X^*$ is the optimal log-monetary reserve in \reff{logX1}  controlled by $\alpha^*$, $K$ is the solution to the BSRE:
\beqs
dK_t &=& \big[ 2(\kappa + q)(I_t) K_t - 2 K_t^2   -   \frac{1}{2} (\eta-q^2)(I_t) \big] dt + Z_t^K dW_t^0, \;\; 0 \leq t \leq T, \; K_T = \frac{c}{2},
\enqs
$\Lambda$ is the solution to the BSRE
\beqs
d\Lambda_t &=&  2 \Lambda_t^2 dt + Z_t^\Lambda dW_t^0, \;\;\; 0 \leq t \leq T, \; \Lambda_T = 0, 
\enqs
and $Y$ is the solution to the linear BSDE
\beqs
dY_t &=& [ 2 \Lambda_t Y_t  - 2 \sigma(I_t) \rho(I_t)  Z_t^Y] dt + Z_t^Y dW_t^0, \;\;\; 0 \leq t \leq T, \; Y_T = 0. 
\enqs
The nonnegative solution $K$ to the BSRE is, in general,  not explicit, while the  solution for $(\Lambda,Y)$ is  obviously equal to 
$\Lambda$ $\equiv$ $0$ $\equiv$ $Y$. From \reff{optiinter3}, it is then clear that  $\E[\alpha_t^*(X_t^*)|W^0]$ $=$ $0$, so that  
the conditional mean of the optimal log-mo\-netary reserve is governed from  \reff{logX1} by
\beqs 
d \E[X_t^*|W^0] &=&   \sigma(I_t)  \rho(I_t) dW_t^0. 
\enqs
The optimal control can then be expressed pathwise  as 
\beqs
\alpha_t^*(X_t^*)  &=& - (2K_t +q(I_t))(X_t^* - x_0 - \int_0^t \sigma(I_s) \rho(I_s)  dW_s^0), \;\;\; 0 \leq t \leq T. 
\enqs

\end{document}